\documentclass[preprint,12pt,authoryear,nopreprintline]{elsarticle}

\usepackage{graphicx}

\usepackage{amssymb}
\usepackage{bm}
\usepackage{siunitx}

\usepackage{lineno}

\usepackage[utf8]{inputenc}
\usepackage{underscore}
\usepackage{todonotes}
\usepackage[colorlinks]{hyperref}
\usepackage{nomencl}
\usepackage[c1]{optidef}
\usepackage{float}
\usepackage{booktabs}
\usepackage{bm}
\usepackage{comment}

\newcommand*\rot{\rotatebox{90}}

\usepackage[version=4]{mhchem}

\newcommand{\R}{\ensuremath{\mathbb{R}}}

\makenomenclature

\usepackage{ifthen}
  \renewcommand{\nomgroup}[1]{%
  \item[\bfseries
  \ifthenelse{\equal{#1}{T}}{Abbreviations and Acronyms}{%
  \ifthenelse{\equal{#1}{S}}{Symbols}{}}%
  ]}

\usepackage{multirow}
\usepackage{array}
\newcolumntype{X}{>{\global\let\currentrowstyle\relax}}
\newcolumntype{D}{>{\currentrowstyle}}

\journal{Computers \& Chemical Engineering}

\begin{document}

\begin{frontmatter}

	\title{Why fixing alpha in the NRTL model might be a bad idea -- Identifiability analysis of a binary Vapor-Liquid equilibrium.}

	\author[1]{Volodymyr Kozachynskyi\corref{cor}}
	\cortext[cor]{Corresponding author}
	\ead{volodymyr.kozachynskyi@tu-berlin.de}
	\author[1]{Christian Hoffmann}
	\ead{c.hoffmann@tu-berlin.de}
	\author[1]{Erik Esche}
	\ead{erik.esche@tu-berlin.de}

	\affiliation[1]{
        organization={Technische Universität Berlin, Process Dynamics and Operations Group},
        addressline = {Stra\ss e des 17.~Juni 135},
		city={Berlin},
		postcode={10623},
		country={Germany}}

	\begin{abstract}
        New vapor-liquid equilibrium (VLE) data are continuously being measured and new parameter values, e.g., for the nonrandom two-liquid (NRTL) model are estimated and published.
        The parameter $\alpha$, the nonrandomness parameter of NRTL, is often heuristically fixed to a value in the range of 0.1 to 0.47.
        This can be seen as a manual application of a (subset selection) regularization method.
        In this work, the identifiability of the NRTL model for describing the VLE is analyzed.
        It is shown that fixing $\alpha$ is not always a good decision and sometimes leads to worse prediction properties of the final parameter estimates.
        Popular regularization techniques are compared and Generalized Orthogonalization is proposed as an alternative to this heuristic.
        In addition, the sequential Optimal Experimental Design and Parameter Estimation (sOED-PE) method is applied to study the influence of the regularization methods on the performance of the sOED-PE loop.
	\end{abstract}


\end{frontmatter}

\printnomenclature

\nomenclature[T]{PE}{Parameter Estimation}%
\nomenclature[T]{OED}{Optimal Experimental Design}%
\nomenclature[T]{sOED-PE}{sequential Optimal Experimental Design and Parameter Estimation}%
\nomenclature[T]{NRTL}{Nonrandom Two-Liquid}%
\nomenclature[T]{VLE}{Vapor-Liquid Equilibria}%
\nomenclature[T]{FFD}{Full Factorial Design}%
\nomenclature[T]{FIM}{Fisher Information Matrix}%
\nomenclature[T]{WLS}{Weighted Least Squares}%
\nomenclature[T]{DOF}{Degree Of Freedom}%
\nomenclature[T]{sv}{Singular Value}%
\nomenclature[T]{SCE}{Scenario}%

\nomenclature[S]{$C$}{parameter variance--covariance matrix}
\nomenclature[S]{$C_{M}$}{measurement error variance--covariance matrix}
\nomenclature[S]{$h$}{set of relations between measured and state variables}
\nomenclature[S]{$S$}{sensitivity matrix}
\nomenclature[S]{$\overline{S}$}{scaled sensitivity matrix}%
\nomenclature[S]{$H_{\theta}$}{Hessian Matrix regarding model parameters}
\nomenclature[S]{$\Phi$}{objective function}

\nomenclature[S]{$\mu$}{experiment index}
\nomenclature[S]{$N_{\mu}$}{number of experiments}
\nomenclature[S]{$u$}{control variables}
\nomenclature[S]{$v$}{eigenvector}
\nomenclature[S]{$N_u$}{number of control variables}
\nomenclature[S]{$N_{\theta}$}{number of parameters}
\nomenclature[S]{$N_y$}{number of measured response variables}
\nomenclature[S]{$\theta$}{model parameters}
\nomenclature[S]{$\hat\theta$}{estimated parameter values}
\nomenclature[S]{$\theta^\ast$}{``true'' parameter values}
\nomenclature[S]{$x$}{state variables}
\nomenclature[S]{$x^{V}$}{vapor mole fraction}
\nomenclature[S]{$x^{L}$}{liquid mole fraction}
\nomenclature[S]{$T$}{temperature}
\nomenclature[S]{$P$}{pressure}
\nomenclature[S]{$y$}{predicted response variables}
\nomenclature[S]{$e$}{residuals}
\nomenclature[S]{$y^m$}{measured response variables}
\nomenclature[S]{$Y^m$}{total experimental data matrix}
\nomenclature[S]{$Y$}{total predicted response matrix}
\nomenclature[S]{$\varepsilon$}{measurement error}
\nomenclature[S]{$\sigma_y^M$}{standard deviation of measurement uncertainty}
\nomenclature[S]{$\sigma_y$}{standard deviation of prediction uncertainty}
\nomenclature[S]{$\Sigma_y$}{measurement error standard deviation matrix}
\nomenclature[S]{$\Sigma^{p}_{p}$}{prediction error standard deviation matrix}
\nomenclature[S]{$f$}{model equations}

\nomenclature[S]{$\varsigma_i$}{singular value}

\nomenclature[S]{$s_y$}{measurement error estimate}
\nomenclature[S]{$\lambda$}{eigenvalue}
\nomenclature[S]{$T^{sp}$}{measurement sampling time grid}
\nomenclature[S]{$T^{u}$}{control switching time grid}
\nomenclature[S]{$\alpha$}{nonrandomness NRTL parameter $\alpha$}
\nomenclature[S]{$\beta$}{confidence level}
\nomenclature[S]{$N_{MC}$}{number of Monte Carlo simulations}

\nomenclature[S]{$U$}{matrix of control values for experimental data}
\nomenclature[S]{$U^{p}$}{matrix of control values for model prediction}

\nomenclature[S]{$Q^{0.95}_{y}$}{prediction quality metric}
\nomenclature[S]{$Q^{0.95}_{\theta}$}{linearity metric}
\nomenclature[S]{$\overline{W}$}{normality metric}
\nomenclature[S]{$D^{U}$}{discrepancy metric}

\nomenclature[S]{$A$}{NRTL parameter or OED criterion}
\nomenclature[S]{$B$}{NRTL parameter}

\section{Introduction}

Accurate prediction of vapor-liquid equilibria (VLE) is critical for the reliability of many complex chemical engineering models and therefore for the design and operation of chemical processes.
Every year, the database of VLE experimental data grows steadily, new mixtures are measured, old measurements are replaced by new ones, and new parameter estimates are published.
The nonrandom two-liquid (NRTL) model~\citep{Renon1968} is one of the most popular models for describing phase equilibria: it is applicable to partially and completely miscible systems and can accurately describe different phase equilibria for multicomponent systems~\citep{Prausnitz1999}.
Studies that measure new equilibrium data typically focus on comparing new data with previously published data, thermodynamic consistency tests, and comparison of different models that were estimated using new data, e.g., NRTL vs.~UNIQUAC.
Estimation of parameters is usually not discussed in detail in the area of thermodynamics as it is considered a well-established procedure.
However, the NRTL model is strongly nonlinear and the model parameters are highly interconnected~\citep{Labarta2022}, making parameter estimation more difficult than simple linear regression.

For a binary mixture, there are at least three parameters to estimate: two binary interaction parameters and the nonrandomness parameter $\alpha$.
The parameter $\alpha$ renders the NRTL model flexible and applicable to many different cases.
However, it is hard to find a publication where the parameter $\alpha$ is estimated and not fixed to some value based on heuristics proposed in the original literature~\citep{Renon1968}.
Both statements sound contradictory: why fix the parameter that makes the model more flexible?

In this work, we analyze a binary VLE model decribed with NRTL using identifiability analysis to gain insight into (i) why the parameter $\alpha$ is usually fixed and what consequences this might have for Parameter Estimation (PE) and Optimal Experimental Design (OED) optimization problems, and (ii) whether regularization methods are a feasible alternative to fixing $\alpha$.
We use these new insights to consider sequential Optimal Experimental Design and Parameter Estimation (sOED-PE) as replacement for the still widely used Full Factorial Design (FFD).

The article continues with a review of current research on the NRTL model and parameter estimation, identifiability analysis, regularization methods, as well as sOED-PE. 
Afterwards, the theory for both PE and OED optimization problems is presented along with four popular regularization techniques. 
In the methodology section, a Monte Carlo-based approach and metrics for comparison are proposed for the analysis of the regularization techniques.
The results are presented in two parts. 
First, the consequences of fixing $\alpha$ are shown based on the PE performed on the typical FFD measurement grid.
Second, the influence of the regularization techniques on the sOED-PE loop is analyzed.

\section{State of the art}

Relevant research on the NRTL model itself, identifiability analysis for parameter estimation, and optimal experimental design is presented in the following.

\subsection{NRTL model}

NRTL, originally published by \cite{Renon1968}, is a local composition model used to describe various types of equilibria for multiple components.
In the simplest case of a binary mixture, it reads as~\citep{Prausnitz1999}:

\begin{equation} \label{eq:nrtl}
\begin{gathered}
\ln \gamma_1=x^2_2\left[\tau_{21}\left(\frac{G_{21}}{x_1+x_2 G_{21}}\right)^2 +\frac{\tau_{12} G_{12}} {(x_2+x_1 G_{12})^2 }\right] \\
\ln \gamma_2=x^2_1\left[\tau_{12}\left(\frac{G_{12}}{x_2+x_1 G_{12}}\right)^2 +\frac{\tau_{21} G_{21}} {(x_1+x_2 G_{21})^2 }\right] \\ 
G_{12}= \exp (-\alpha_{12}\ \tau_{12}) \quad G_{21}=\exp (-\alpha_{21}\ \tau_{21}) \\
\end{gathered}
\end{equation}

Binary interaction parameters $\tau$ are usually temperature dependent, while $\alpha_{12}$ is equal to $\alpha_{21}$, henceforth referred to as $\alpha$.
In the following, $\tau$ is assumed to be linearly dependent on the temperature.
The formulation in~\autoref{eq:tau} is usually extended with additional parameters to describe liquid-liquid equilibria.

\begin{equation} \label{eq:tau}
    \tau_{12} = A_{12} + \frac{B_{12}}{T} \quad \tau_{21} = A_{21} + \frac{B_{21}}{T} \\
\end{equation}

For a binary mixture, at least five parameters must be estimated: $A_{12}, B_{12}, A_{21}, B_{21}$ and $\alpha$.
This is in itself a complex optimization problem because (i) multiple, non-unique parameter sets may exist that are equally suitable to describe the VLE~\citep{Gau2000}, and (ii) the parameters are highly correlated, as seen in Equations~\ref{eq:nrtl} and~\ref{eq:tau}.
The first problem may be solved by using global optimization methods~\citep{McDonald1995,Bollas2009}, or by determining the multiple feasible areas that lead to optimal solutions~\citep{Werner2023}.
Although global optimization methods play an important role in finding the best parameters, they are not the focus of this work because we want to analyze the performance of local optimization solvers that are initialized close to the global solution.
The second problem, the interdependence of the NRTL parameters, is directly related to nonlinearity and negatively affects the convergence of local optimization methods, which is a typical problem for nonlinear regression~\citep{Ratkowsky1990}.
The influence of nonlinearity can be reduced by higher accuracy of measurement data, reparametrization, a very good initial guess for the estimated parameters, the use of regularization methods, or the selection of experimental conditions maximizing information content.

Several authors have already proposed to rewrite the original NRTL model to reduce the number of estimated parameters, e.g., \cite{Gebreyohannes2014}. 
Despite progress in this area, we often see that additional thermodynamic models are only added to process simulators when there is sufficient demand for them.
Others have analyzed the combination of different parameter values to obtain a good initial guess for the parameter estimation \citep{Serafimov2002, Labarta2022}.
Such an analysis of the NRTL model is important to gain a better understanding of the NRTL model and its internal mechanisms, but it only helps to find the parameter estimates, similar to global optimization methods, and thus is not the focus of this work.

The technically simplest and most common solution is to fix parameter $\alpha$.
It is usually done based on the rules proposed by \citet{Renon1968}.
In recent years, new research has been published that provides more insight into the influence of the values of $\alpha$ \citep{Klerk2023} and new heuristics for selecting the fixed value.
In practice, this method can be categorized as a manual application of the subset selection-based regularization algorithm.
However, subset selection is only one part of the family of regularization algorithms available for parameter estimation.
Despite the existence of regularization algorithms for several decades, they are rarely used in thermodynamics, e.g.,~\citep{Hoffmann2019}. 
In the following, we briefly discuss identifiability analysis and available regularization techniques for parameter estimation.

\subsection{Ill-posed problems and regularization}

Parameter estimation problems are often ill-posed inverse problems \cite{Lopez2015}.
Ill-posedness is usually understood as a violation of one of the three conditions for well-posed problems and was formalized later by \cite{Kabanikhin2008}. These can be summarized as (i) existence of a solution, (ii) uniqueness of a solution, and (iii) continuity of the solution space.
Parameter estimation problems are usually ill-posed because of the ill-conditioning, i.e., the system matrices are ill-conditioned~\citep{Kabanikhin2008}.
The ill-posedness of parameter estimation problems leads to the non-identifiability of parameters, i.e., no values or no unique values for the model parameters can be determined based on the available experimental data.

In the past, several groups have suggested methods to improve the conditioning of parameter estimation problems, or to transform an ill-posed problem into a well-posed problem, to ensure that at least some of the model's parameters can at least be properly estimated.
Examples of these methods are Subset Selection, Tikhonov regularization or Truncated Singular Value Decomposition (TSVD).

\subsubsection{Subset Selection}

In Subset Selection (SsS), the parameters are partitioned into well-conditioned (or identifiable) and ill-conditioned (or unidentifiable) subsets~\citep{Burth1999, VelezReyes1995, Fink2007}.
There are several different implementations of SsS, but they usually share a number of recurring elements:
The sensitivity matrix $S$ is computed, which is the derivative of the state variables relating to measurements with respect to the model parameters.
Then a rank revealing decomposition, e.g., Cholesky decomposition, is applied to determine the rank of $S$.
In case $\text{rank}\left(S\right)$ equals $m$, where $m$ is a number of parameters, all parameters are identifiable.
For a rank below $m$, the linearly independent columns of $S$ are grouped into a reduced sensitivity matrix, for which the condition number $\kappa$ and the spread of the parameters are below a predefined thresholds.
The parameters are organized accordingly leading to a split between identifiable parameters (corresponding to the linearly independent columns of $S$) and the unidentifiable ones.
This regularization is performed repeatedly in combination with successive parameter estimation.
The parameter estimation problems that operate only on these reduced sets of identifiable parameters should therefore be well-conditioned.
Variations in the implementations appear in the computation of the sensitivity matrix, the choice of rank revealing decomposition, etc.

Subset selection has found numerous applications for the regularization of parameter estimation problems by now, ranging from HIV models~\citep{Fink2007}, across electric machinery~\citep{Burth1999, VelezReyes1995}, an activated sludge model~\citep{Kim2010}, and various other case studies of biological and biochemical systems~\citep{Lopez2015}, as well as applications in chemical engineering~\citep{Mueller2014}.\\

\subsubsection{Tikhonov regularization}

In Tikhonov regularization, \textit{a priori} information about the parameter values is added to the objective function of the parameter estimation problem.
If a likely value for a particular parameter is known, a penalty term is added to the objective function to induce a bias towards that value~\citep{Hansen1998}.
Various formulations exist here as well, a popular one is given in \autoref{eq:TikhPenalty}, wherein $\theta$ are the parameters to be estimated, $\theta^R$ are prior known values, $L$ is a diagonal matrix, which may contain the confidence in the priors (the entry is set to 0 to enforce no penalty), and $\lambda$ is the overall positive scalar parameter to control the strength of the regularization.
The penalty term $f^\textrm{penal}(\theta)$ in \autoref{eq:TikhPenalty} is then added to the objective function of the PE problem.
\begin{gather} \label{eq:TikhPenalty}
    \begin{aligned} 
    f^\textrm{penal}(\theta) &:= \frac{\lambda^2}{2} \cdot \Omega(\theta)\\
    \Omega(\theta) &:= \left\| L \cdot (\theta - \theta^R) \right\|_2^2
\end{aligned}
\end{gather}
The larger $\lambda$ is chosen, the more biased the solution is towards $\theta^R$, but the optimization problem is also smoother and better conditioned at the same time.
Variations in the implementation exist through different functions for $\Omega(\theta)$ and different choices regarding the appropriate selection of, e.g., $\lambda$ and $L$.
Examples of Tikhonov regularization applied to parameter estimation problems can be found for models of geothermal heat exchangers~\citep{Duff2017} and also for liquid chromatography~\citep{Barz2016}.\\

\subsubsection{Truncated singular value decomposition}

If TSVD is used as regularization of the parameter estimation problem, the objective function and the search space of the parameters remain unchanged. TSVD is applied on $S$ and a well-conditioned, but rank-deficient version of $S$ is obtained by setting small non-zero singular values below a chosen threshold to zero. The resulting $S^{reg}$ is used to compute Jacobian and Hessian matrices to solve the parameter estimation problem.
TSVD is also used in~\citep{Lopez2015} for biological or bioprocessing examples and in~\citep{ELSheikh2013} for subsurface flow models.\\

\subsection{Optimal Experimental Design}

As discussed earlier, it is possible to improve the identifiability of the model parameters by conducting experiments that maximize information content regarding the uncertain model parameters.
Almost all available binary VLE data have one thing in common -- experimental conditions are selected based on FFD.
There are good reasons for performing experiments this way, e.g., it is easier to apply thermodynamic consistency tests as they require the (numerical) solution of an integral.
If we focus on the parameter estimation only, however, a viable alternative to FFD is model-based OED.

In classical OED formulations, e.g., A, D, E criteria, the objective is usually a function of the parameter sensitivity matrix $S$~\citep{Franceschini2008}.
New OED objective functions and modifications of the existing OED methods are continuously being suggested,~\citep{Fleitmann2022, Cenci2023} but none of them, even the standard ones, are popular in the design of experiments for equilibrium data.
Already ten years ago, it was shown by~\citet{Dechambre2014} for liquid-liquid equilibrium measurements that model-based OED techniques are applicable to the thermodynamic measurements.
Later,~\citet{Duarte2021} showed the applicability of OED for VLE measurements using artificial data.
This was confirmed by~\citep{Bubel2024}, who showed by comparing FFD and sOED-PE on real experimental data that sOED-PE needed only half the number of experiments as FFD to achieve similar prediction accuracy and quality.

Despite its advantages, OED is highly dependent on the conditioning of the sensitivity matrix, i.e., the value of the currently known parameter values and their uncertainty.
The closer the current parameter estimate is to the unknown, ``true'' values of the parameters, the better the performance of OED methods.
This undesirable property of OED is well known~\citep{Franceschini2008}, and is usually handled by combining OED with PE in a sequential manner~\citep{Kozachynskyi2024}.
After initial experiments are performed, initial parameter estimates are obtained.
These estimates are used to solve an OED problem and to suggest new experiments.
The sequential execution of PE and OED problems is further referred to as sOED-PE iterations.
In general, sOED-PE iterations are repeated until certain criteria, such as target parameter or prediction accuracies, are met.
It should be noted that an ill-conditioned sensitivity matrix affects the performance of not only PE but also OED problems, especially during the first sOED-PE iterations~\citep{Lopez2015}, so the use of regularization algorithms is important for the convergence of sOED-PE.

OED based on sensitivity matrix $S$, SsS, and TSVD regularization techniques operate under important assumptions that will be described in detail later: the parameter uncertainty is approximated by the sensitivity matrix and is equal to the lower bound of Cram\'er-Rao.
This means that the parameter estimate is always assumed to be fully efficient and unbiased, i.e., the model is linear in parameters, there is no plant-model mismatch, and the experimental data are not biased.

\section{Theory}

In the following, we focus on the use of steady-state models formulated as:

\begin{equation} \label{eq:model}
	\begin{gathered}
		0 = f(x, u, \theta) \\
	\end{gathered}
\end{equation}

$x \in \R^{N_{x}}$ is the vector of state variables, $u \in \R^{N_{u}}$ is the vector of control variables, $\theta \in \R^{N_{\theta}}$ is the vector of model parameters.

The elements of the vector of observable predicted response variables $y \in \R^{N_{y}}$ depend on the solution $x$ of the model shown in Eq. \ref{eq:model} and the measurement function $h$:

\begin{equation}
	y = h(x)
\end{equation}

Experimental data consists of two matrices: experimental controls $U = \lbrack u_1, \ldots, u_{N_{\mu}} \rbrack ^T \in \R ^ {N_{\mu} \times N_{u}}$ and corresponding measurement data $Y^{m} = \lbrack y^{m}_{1}, \ldots, y^{m}_{N_{\mu}} \rbrack ^T \in \R ^ {N_{\mu} \times N_{y}}$, where $N_{\mu}$ is the number of different experiments.
The model (Eq. \ref{eq:model}) must be solved for each experimental control $u_{\mu} \in U$ yielding a corresponding matrix of model responses 
$Y = \lbrack y(u_{1}, \theta), \ldots, y(u_{\mu}, \theta) \rbrack ^T \in \R ^ {N_{y} \times N_{\mu}}$.

The difference between the model prediction and the measurement data is called the residual $e_{\mu,i}$:

\begin{equation}
	\begin{gathered}
        e_{\mu,i} = y_{\mu,i} - y^{m}_{\mu,i} \\
	\end{gathered}
\end{equation}

Assuming structural correctness of the model and a large number of measurements, model evaluation at the ``true'' parameter values $\theta^\ast$ leads to the residual of each measurement being normally distributed around zero with the standard deviation of the measurement noise $\varepsilon_{\mu,i}$.
The sensitivity of the predicted response variables to the parameters is collected in the sensitivity matrix $S(U, \theta) = dY(U,\theta) / d\theta \in \R^{N_{y}\cdot N_{\mu}\times N_{\theta}}$.

\subsection{Parameter Estimation}

Assuming that the standard deviation of the measurement noise $\varepsilon_{\mu,i}$ is the same for all experiments, the objective function of the PE can be written in the Weighted Least Squares form as \citep[p.56]{Bard1974}:

\begin{equation} \label{eq:wlsq}
    \Phi^{\text{WLS}} = \sum_{\mu=1}^{N_{\mu}}e_{\mu}^T  C_{M}^{-1} e_{\mu}
\end{equation}

where the weighting matrix 
$C_{M}^{-1} \in \R^{N_{y} \times N_{y}}$
is the inverse of the measurement variance-covariance matrix $C_{M}$.
It is assumed that measurement noise is random, drawn from a normal distribution with zero mean, and is independent for different measurements and experiments. 
Thus, all non-diagonal elements of the $C_{M}$ matrix are zero, and the diagonal elements are the variances of each measurement ${\sigma^M_y}^2$.
The square root of the diagonal is the standard deviation matrix $\Sigma_{y}$, i.e., $\Sigma_{y} = \sqrt{\operatorname{diag}\left(C_{M}\right)}$.

Given the above assumptions, the optimal solution $\hat{\theta}$ is the maximum likelihood estimate and is asymptotically unbiased \citep[p.61]{Bates1988}.
Noise in the measurements leads to parameter uncertainty, which according to Cram\'er-Rao's bound is greater or equal to the inverse of Fisher information matrix (FIM), estimated by the parameter variance--covariance matrix $C(U, \hat{\theta}) \in \R^{ N_{\theta}\times N_{\theta}}$.
Assuming that $\hat\theta$ is close to $\theta^\ast$, the matrix $C(U, \hat{\theta})$ can be approximated by the inverse of the Hessian of the objective function $H_{\theta}=\nabla_\theta^{2} \Phi^{\text{WLS}}$. 
Assuming linearity of the model, the Hessian matrix can be computed using the sensitivity matrix according to Gauss' approximation, so that:

\begin{equation} \label{eq:parcov}
    C(U, \hat\theta) \geq H_{\theta}^{-1}(U, \hat\theta) \approx \left(S(U, \hat\theta)^T C_{M}^{-1} S(U, \hat\theta) \right)^{-1}
\end{equation}

The confidence in the model predictions, the prediction accuracy, is described by the standard deviation of the prediction $\sigma_{y}$, which for given controls $u$ under linearization assumptions can be approximated as~\citep{Bates1988}:

\begin{equation} \label{eq:prediction-std}
	\begin{gathered}
        \sigma_y = \operatorname{diag}\left(\Sigma_y \sqrt{S(u,\hat\theta)^T C(U,\hat\theta) S(u,\hat\theta)}\right) \\ 
        y(u, \hat\theta) \pm \sigma_{y} \cdot t(\text{DOF}, \beta / 2)
	\end{gathered}
\end{equation}

$t(\text{DOF}, \beta / 2)$ is the upper $\beta / 2$ quantile for Student's distribution with degrees of freedom calculated as \citep{Bard1974}:

\begin{equation}\label{eq:dof}
    \text{DOF} = N_{\mu} - N_{\theta} / N_{y} 
\end{equation}

Thus, parameter accuracy, the confidence in the parameter estimates described by the parameter variance--covariance matrix, directly affects the prediction accuracy of the model.

\subsection{Optimal Experimental Design}

As seen in \autoref{eq:prediction-std}, there are two ways to influence model prediction accuracy: either by improving measurement accuracy $\Sigma_y$, or by selecting experimental controls $U$ that lead to a ``smaller'' overall parameter variance-covariance matrix $C_(U, \hat\theta)$.
The latter can be achieved by solving the OED optimization problem:

\begin{equation} \label{eq:OED}
    \begin{gathered}
        \min_{U^{new}}~\Phi^{crit} \left( S\left(U, \theta\right) \right) 
    \end{gathered}
\end{equation}
where $U$ is the matrix of controls of experiments which are already conducted and proposed for execution by OED: $U = U^{exp} \cup U^{new}$.
If $U^{new}$ has more than one row, multiple experiments are computed simultaneously.

Several OED criteria exist to formulate $\Phi^{crit}$ with the goal of improving parameter precision \cite{Franceschini2008}.
Criteria A, D, and E are the most commonly used, where objective $\Phi^A$ minimizes the trace of $C_{U, \theta}$; objective $\Phi^D$ minimizes the determinant of $C_{U, \theta}$; objective $\Phi^E$ minimizes the largest eigenvalue of $C_{U, \theta}$~\citep{Kozachynskyi2024}. 

\subsection{Regularization}

As discussed above, identifiability, is an issue for PE and indirectly also for OED optimization problems.
Subset selection regularization algorithms are the simplest type to implement, and compared to other methods, parameter estimates and their uncertainty can be computed without modification.
The idea is to fix the parameters that are selected by the algorithm as non-identifiable, and not to estimate them.

Regularization methods use scaled parameter sensitivity matrices, in order to eliminate the difference in scale between the measurement and parameter variables.
There are several approaches to scaling the sensitivity matrix, only two of which are used in this work.
In both cases, the matrix is scaled by the parameter values, to adjust the sensitivity of the parameter to its absolute value.
To adjust the magnitude of the measurement values, \cite{Yao2003} suggests using the evaluated values of the measurement variables $y_{\mu,i}$, while other methods use scaling based on the measurement noise $\Sigma_{y}$:

\begin{gather}
    \overline{S}_{Yao} = S \cdot \text{diag}[\hat\theta] \cdot y^{-1} \\
    \overline{S} = S \cdot \text{diag}[\hat\theta] \cdot \Sigma_y^{-1}
\end{gather}

Four popular regularization algorithms are briefly introduced below.
For a more detailed description of these algorithms, please refer to the original works.

\subsubsection{Eigenvalue algorithm (E)}

In their study,~\cite{Quaiser2009} propose an approach for analyzing the identifiability of parameters by examining the eigenvalues of the Fisher information matrix. 

\begin{enumerate}
    \item Compute the FIM matrix, denoted by $\text{FIM} = S^{T} S$.
    \item Compute the eigenvalues $\lambda^{j}$ and the eigenvectors $v^{j}$ of the FIM matrix, sorted from the smallest to the largest eigenvalues.
    \item If the smallest eigenvalue, $\lambda^{1}$, is greater than or equal to a threshold, $\varepsilon=0.001$, all parameters are identifiable.
    \item Otherwise, the parameter corresponding to the highest absolute value in the eigenvector $v^{1}$ is marked as non-identifiable and step 1 is repeated for all other parameters.
\end{enumerate}

\subsubsection{SVD algorithm (SVD)}

In their publication,~\citet{Lopez2015} suggests executing the identifiability analysis directly using a scaled sensitivity matrix. 
In consecutive work,~\cite{Lopez2015} suggest different scaling techniques for the sensitivity matrix while maintaining the same algorithm:

\begin{enumerate}
    \item Conduct singular value decomposition on the scaled sensitivity matrix $\overline{S}$, sorting the singular values (sv) $\varsigma_i$ from the largest to the smallest value.
    \item Compute the condition number $\kappa$ as the ratio of the largest and smallest SVs.
    \item If the condition number is less than a threshold $\varepsilon = 1000$, all parameters are identifiable.
    \item Otherwise, identify the largest number $r$ of identifiable parameters that result in a  the condition number less than $\varepsilon$.
    \item Use QR decomposition with column pivoting to obtain the identifiability order of the parameters, select $r$ parameters as identifiable ones.
\end{enumerate}

\subsubsection{Forward selection algorithm (FS)}

\cite{Yao2003} suggested the use of the magnitude of the sensitivity matrix with an orthogonalization method to obtain a set of identifiable parameters.

\begin{enumerate}
    \item Compute the magnitude of each column of the scaled sensitivity matrix $\overline{S}_{Yao}$, to determine the most identifiable parameter.
    \item Use Gram–Schmidt orthogonalization to find and remove the most identifiable parameter from the sensitivity matrix.
    \item Repeat step two with the modified sensitivity matrix until the largest magnitude of the reduced sensitivity matrix is less than the threshold $\varepsilon = 0.04$. 
\end{enumerate}

\subsubsection{Generalized orthogonalization algorithm (GO)}

In their work,~\cite{Chu2012} propose to select identifiable parameters by solving an optimization problem.

\begin{enumerate}
    \item Calculate the scaled sensitivity matrix $\overline{S}$ of the currently analyzed parameter estimate vector $\hat\theta$.
    \item Use a brute-force search algorithm to find a subset of parameters $J$ that corresponds to the biggest determinant of the FIM matrix given by:
    $\operatorname{det}~\left(\overline{S}\left(\colon,J\right)^{T} \overline{S}\left(\colon,J\right)\right)$
\end{enumerate}

\subsubsection{Overview}

The selected algorithms are among the most cited algorithms used in process engineering and use various criteria for regularization. 
The \emph{FS} algorithm uses Gram-Schmidt orthogonalization for parameter ranking, which is also used in a more general way in the \emph{GO} algorithm.
The \emph{E} algorithm ranks parameters based on the eigenvalues of the FIM matrix and selects the identifiable ones based on the eigenvalue threshold, while the \emph{SVD} algorithm ranks parameters based on the singular value of the sensitivity matrix and uses the condition number as the parameter identifiability threshold.

Apparently, the outcomes of the three methods depend on the user-defined thresholds.
The \emph{SVD} algorithm uses a threshold for maximum condition number that can be mathematically justified.
There is, however, no apparent methodology for selecting the thresholds of the \emph{E} and \emph{FS} algorithms.
The GO algorithm is the only algorithm that does not require the specification of a threshold, but it becomes impractical with a large number of parameters.

\section{Methodology}\label{methodology}

In order to analyze the influence of fixing $\alpha$ and different regularization techniques, and to evaluate the applicability of the linearization assumptions introduced in Eqs.~\ref{eq:parcov} and~\ref{eq:prediction-std}, it is important that the measurement data are not biased and that the ``true'' parameter values are known.
The Monte Carlo method~\citep[p.46]{Bard1974} ensures these properties and is usually used to quantify the ``exact'' parameter and prediction uncertainties and the performance of sOED-PE methods~\citep{Kozachynskyi2024}.
Usually Monte Carlo is used without regularization techniques, so in the following the original Monte Carlo method is presented together with the modified versions and the metrics that are further used to compare the influence of different regularization techniques on PE, OED and sOED-PE.

\subsection{Monte Carlo method \emph{$\left(\text{MC}\right)$}}\label{MC}

The ``exact'' parameter confidence region is generated as follows:

\begin{enumerate}
    \item The ``true'' parameter values $\theta^\ast$ are selected.
    \item The measurement grid $U$ is selected for the generation of artificial data.
    \item The prediction grid $U^{p}$ is selected for generation of the ``true'' data $Y^{p}_{true}$.
    \item The ``true'' model is evaluated at measurement grid $U$, resulting in the generation of  ``true'' measurement data for selected measurement variables $y_i$.
    \item The ``true'' data is perturbed using noise drawn from normal distribution for a given measurement variance--covariance matrix $C^{m}$.
    \item The parameter estimate $\theta_{MC}$ is calculated by minimizing the $\Phi^{\text{WLS}}$ (Eq. \ref{eq:wlsq}).
\item The measurement error estimate vector is calculated as $s_{y} = \left(\sqrt{\text{diag}\left({e_{\mu,i}}^Te_{\mu,i}\right)} \cdot \text{DOF}\right) \in \R^{N_{i}}$ to get an estimate of the measurement error standard deviation $\Sigma_y$. 
    \item The model is evaluated at the prediction grid $U^{p}$ with the parameter estimate $\theta_{MC}$ to obtain the model prediction $Y^{p}$.
    \item Steps 4 till 6 are repeated $N_{MC}$ times, resulting in a final matrix of the parameter estimates $\bm{\theta} \in \R^{N_{\theta}\times N_{MC}}$, measurement error estimates $\bm{s_{y}} \in \R^{N_{i} \times R^{N_{MC}}}$tensor of the model predictions $\bm{Y^{p}} \in \R^{N_{\mu} \times N_{y} \times N_{MC}}$ with corresponding standard deviations  $\bm{ \Sigma^{p}_{y}} \in \R^{N_{\mu} \times N_{y} \times N_{MC}}$ (Eq. \ref{eq:prediction-std}).
\end{enumerate}

The resulting matrix of parameter estimates and tensor of the model predictions can now be used to test whether the linearization assumption holds.
The matrix of control values used for model prediction $U^{p}$ can be different from the measurement grid $U$ to check the validity of the parameters in the design space where no experimental data were used for estimation.
According to the definition, only a subset of the parameters in $\bm{\theta}$, defined by probability $\beta$, should contain true parameter values $\theta^{\ast}$ within their respective linearized confidence intervals~\citep{Bates1988}:

\begin{equation} \label{eq:confidence-interval-parameter}
    \theta \pm \sigma_{\theta} = \theta \pm \Sigma_y \sqrt{\text{diag}~C(U,\theta)} \cdot \left(\text{DOF}, \beta / 2\right)
\end{equation}

Similarly, not all parameters in $\bm{\theta}$ result in accurate predictions, meaning that ``true'' prediction $Y^{p}_{true}$ is inside the confidence interval $\bm{Y^{P}} \pm \bm{\Sigma^{p}_{y}} \cdot t(\text{DOF}, \beta / 2)$ defined by probability $\beta$.

The ratio of parameters $\bm{\theta}$ located inside the 95 \% confidence interval of the parameter $\theta^\ast$ to the total number of Monte Carlo simulations $N_{MC}$ is referred to as the \emph{linearity} $Q^{0.95}_{\theta}$.
If the model is linear in parameters, the $\sigma_{\theta_{MC}}$ (\autoref{eq:confidence-interval-parameter}) is equal to $\sigma_{\theta^\ast}$, meaning that for a sufficiently large number of Monte Carlo simulations, the linearity metric $Q^{0.95}_{\theta}$ should be close to 0.95.

The ratio of parameters $\bm{\theta}$ with successful 95 \% confidence interval predictions to the total number of Monte Carlo simulations $N_{MC}$ is further referred to as the \emph{prediction quality} $Q^{0.95}_{y}$.
If the linearization assumption is valid and the number of Monte Carlo simulations $N_{MC}$ is sufficiently large, the values of $Q^{0.95}_{y}$ should be close to 0.95.

The tensor of model predictions $\bm{s}_{i}$ indicates the bias in the parameter estimates.
The columns of the matrix should be normally distributed around the mean close to the respective measurement error $\sigma_{y}^M$.
In the event that this is not the case, overfitting or plant-model mismatch may be observed.

An additional metric, further referred to as \emph{normality}, is used to assess the \emph{linearity} of the model.
The columns of the matrix are evaluated by the Shapiro-Wilk test \citep{Shapiro1965, Virtanen2020}, and the mean value $\overline{W}$ across all parameter values is computed.
A value of $\overline{W}$ close to one suggests that the parameters are normally distributed, indicating that the model is linear in the parameters.

\subsubsection{MC for identifiability \emph{$\left(\text{MC}^{\text{I}}\right)$}}\label{MCi}

To analyze the influence of individual parameters on the PE, one of the parameters can be fixed and excluded from the PE, which will affect the final parameter uncertainty, model prediction accuracy, and model bias.
The Monte Carlo method~$\left(\text{MC}\right)$ in ~\autoref{MC} is then modified as follows: prior to parameter estimation ($\text{MC}$, Step 4), one of the model parameters is fixed to either (i) a ``correct'' value $\hat\theta$, or (ii) a ``wrong'' value.
In scenario  (i), one of the parameters is fixed to a ``true'' value, and four remaining parameters are estimated.
If the model is more linear in the remaining four parameters, as evidenced by higher values of the normality $\overline{W}$ and linearity $Q^{0.95}_{\theta}$ metrics, it is likely that the fixed parameter caused the nonlinearity.
In scenario (ii), one of the parameters is fixed at a value different from the ``true'' value, and four remaining parameters are estimated.
If the ``wrong'' parameter set, as judged by $\bm{s_{y}}$, still provides good parameter accuracy, it indicates that the fixed parameter has little influence, i.e., the model is less sensitive to that parameter.

\subsubsection{MC for regularization \emph{$\left(\text{MC}^{\text{Reg}}\right)$}}\label{MCreg}

To compare regularization algorithms, the Monte Carlo method~$\left(\text{MC}\right)$~\ref{MC} is further modified as follows.
Following the parameter estimation ($\text{MC}$, Step 5), an identifiability analysis is applied on the parameter estimate $\theta_{MC}$.
The next steps are than executed only with those parameters that have been selected as identifiable.
Consequently, the bias introduced by the regularization algorithm can be quantified and compared.
It is anticipated that the regularization algorithm will result in a reduction in the mean prediction error $\overline \sigma_{y}$, at the expense of reduced prediction quality $Q^{0.95}_{y}$.

\subsection{MC with sOED-PE \emph{$\left(\text{MC}_{sOED-PE}\right)$}} \label{MCsoed}

To validate the performance of the sOED-PE method, the following modification of the Monte Carlo method is used:

\begin{enumerate}
    \item The prediction grid $U^{p}$ is selected for generation of the ``true'' data $Y^{p}_{true}$.
    \item Measurements $Y$ for the initial experiments $U^{init}$ are generated with the best parameter estimate $\hat\theta$, which is assumed to be ``true''.
    \item ``True'' data is perturbed using noise drawn from a normal distribution for a given measurement variance--covariance matrix $C^{M}$.
    \item Parameter estimate $\theta_{sOED,MC}$ is calculated by minimizing the $\Phi^{\text{WLS}}$ for $Y$ (Eq. \ref{eq:wlsq}).
    \item All metrics from the Monte Carlo method are calculated (Steps 7 and 8 in $\text{MC},~\autoref{MC}$).
    \item An experimental design $u^{new}_{sOED}$ is calculated by OED (Eq.~\ref{eq:OED}).
    \item $Y$ is updated with new measurement data generated and perturbed at $u^{new}_{sOED}$.
    \item Steps 4 through 7 are repeated, until all $N_{sOED}$ iterations are executed.
    \item Steps 1 through 8 are repeated $N_{MC}$ times to obtain values for the mean and standard deviation of the metric associated with each sOED iteration over the entirety of the $N_{MC}$ simulations.
\end{enumerate}

The resulting tensor of all experimental designs $\bm{U}$ can then be analyzed using the central discrepancy metric $D^{U}$~\citep{Zhou2013}.
Given an experimental design tensor (i), in which the designs are clustered around a few regions, and a design tensor (ii), in which the designs are spread evenly over the design space, the value of the discrepancy is higher for the less spread design (i).

\subsubsection{MC with sOED-PE for identifiability \emph{$\left(\text{MC}_{sOED-PE}^{\text{I}}\right)$}}\label{MCsoedi}

To analyze the influence of individual parameters on the sOED-PE, one of the parameters can be fixed and excluded from the PE and OED, which will affect the final parameter uncertainty, model prediction, model bias, and proposed design of experiments.
The Monte Carlo method with sOED-PE~$\left(\text{MC}_{sOED-PE}\right)$ in~\autoref{MCsoed} is then modified as follows: prior to parameter estimation ($\text{MC}_{sOED-PE}$, Step 4), one of the model parameters is fixed to either (i) a ``correct'' value $\hat\theta$, or (ii) a ``wrong'' value, as suggested in $\left(\text{MC}^{\text{I}}\right)$ ~\autoref{MCi}.
It also means, that during the OED ($\text{MC}_{sOED-PE}$, Step 6), the sensitivity of the fixed parameter is not used to calculate the objective function.
As the number of sOED-PE iterations increases, it is expected that the bias caused by fixing the parameter to the ``wrong'' value, scenario (ii), indicated by the low prediction quality $Q^{0.95}_y$, should become more apparent.
The decrease in prediction quality is expected to be greater the more sensitive the model is to the parameter.
In both scenarios, the lower variability of the proposed designs $u^{new}_{sOED}$ is expected because the sensitivity of only four parameters is considered.

\subsubsection{MC with sOED-PE for regularization \emph{$\left(\text{MC}_{sOED-PE}^{\text{Reg}}\right)$}}\label{MCsoedreg}

Similar to $\text{MC}_{sOED-PE}^{\text{I}}$, to analyze the influence of regularization techniques, parameters marked as non-identifiable are fixed before PE ($\text{MC}_{sOED-PE}$, Step 4).
It is expected that during the first sOED-PE iterations with less experimental data, the parameter estimates may change a lot, leading to a different set of identifiable parameters.
However, the smaller the number of identifiable parameters is, the more biased the sOED-PE should become, resulting in a decrease in prediction accuracy over the increasing number of sOED-PE iterations.
To see the influence of the regularization on the OED ($\text{MC}_{sOED-PE}$, Step 6), the non identifiable parameters can be either included (i) or excluded (ii) from the objective function of the OED (Step 6).
If the set of identifiable parameters changes over iterations, more variability in experimental design is expected.

\subsection{Metrics}

The following is a summary of the metrics that were used in the analysis of the results:

\begin{itemize}
    \item Parameter \emph{normality} $\overline{W}$: smaller values suggest that parameters are not normally distributed, thus the model is more nonlinear.
    \item Model \emph{linearity} $Q^{0.95}_{\theta}$: values smaller than 0.95 indicate that the linearity assumptions do not hold. This is analogous to the concept of $\overline{W}$.
    \item \emph{Prediction quality} $Q^{0.95}_{y}$: values smaller than 0.95 suggest that model predictions are not accurate.
    \item \emph{Estimation accuracy} $\overline{s}_{i}$: values that are not close to the real measurement error $\sigma_{y}^M$ indicate that the parameter estimates are biased. The metric is the average of the measurement error estimator $\bm{s_y}$ over all Monte Carlo simulations. 
    \item \emph{Prediction accuracy} $\overline \sigma_{y}$: smaller values indicate more accurate predictions. The metric is the average of the prediction error $\bm{\Sigma^{p}_{y}}$ over all Monte Carlo simulations and all points on the prediction grid $U^{p}$.
    \item \emph{Spread of experiments} $D^{U}$: smaller values indicate that the experimental design is more spread over the design space. This metric is only used for results generated using $\text{MC}_{sOED-PE}$~(\ref{MCsoed}) methods.
\end{itemize}

\section{Case studies}

Two case studies are used for further analysis: VLE data are generated using an FFD, and NRTL parameters are estimated once (\emph{Case Study I}); VLE data are generated using the sOED-PE method, NRTL parameters are reestimated after each experiment (\emph{Case Study II}).
\emph{Case Study I} is used to analyze parameter uncertainty, model accuracy, the influence of fixing $\alpha$, and the application of regularization techniques.
\emph{Case Study II} provides insight into the expected performance of the sOED-PE method for different mixtures -- with and without the use of regularization techniques.
All research data and the software are available online~\citep{Kozachynskyi2024a}.

\subsection{Implementation}

In~\autoref{tab:mixtures}, 12 selected component systems with different azeotropic behavior are shown.
The binary VLE is modeled under the simplified isofugacity condition with an ideal vapor phase.
The fugacities are described using a $\gamma - \phi$ approach~\cite[p.~440]{Tassios1993}, neglecting the influence of the Poynting effect and the fugacity coefficients of the individual components:

\begin{equation}
\begin{gathered}
	x^{V}_{i} \cdot P  =  x^{L}_{i} \cdot \gamma_{i} \cdot P^{s}_{o,i}  \qquad i = \{1, 2\} \\ 
    \gamma_{i} = \operatorname{f}_{\text{NRTL}} \left( x^{L}_i, T, A_{12}, B_{12}, A_{21}, B_{21}, \alpha \right) \\
    x^{V}_{1} + x^{V}_{2} = 1 \qquad x^{L}_{1} + x^{L}_{2} = 1 \qquad P^{s}_{o,i} = \operatorname{f}_{\text{Wagner25}}(T)
\end{gathered}
\end{equation}

where $x^V_i$ and $x^L_i$ are the vapor and liquid mole fractions of the components respectively, $P^{s}_{o,i}$ is the vapor pressure of the pure components described by the NIST Wagner 25 equation, and $\operatorname{f}_{\text{NRTL}}$ is the NRTL model described in Eqs.~\ref{eq:nrtl} and~\ref{eq:tau}.
More detailed information about selected mixtures and source of NRTL parameters is given in~\autoref{tab:mixtures}.
The NRTL parameters are assumed to be ``true'' even though they are only the current best parameter estimate based on available experimental data.

The artificial experiments are conducted under isobaric conditions.
Liquid mole fraction and pressure $[x^L_1, P]$ are control variables $u$, bounded as follows: $x^L_1$ in $[0.01,~0.99]~\si{\mole\per\mole}$, $P$ in $[0.5,~1.5]~\si{bara}$.
Measurement variables $y$ are either the vapor mole fraction $[x^V_1]$ only, or the vapor mole fraction and the equilibrium temperature $\left[x^V_1, T\right]$.
Although temperature is almost always measured together with the molar fraction data, both measurement variables scenarios $\left(\text{SCE}_{y}\right)$ are used to see the influence of the selected variables on PE and OED.
The standard deviation of the measurement noise $[\sigma^M_{x^V_1},~\sigma^M_{T}]$ is either called \emph{precise} $[0.0002~\si{\mole\per\mole},~0.01~\si{K}]$ or \emph{default} $[0.001~\si{\mole\per\mole},~0.03~\si{K}]$.
Both scenarios are referred to as measurement accuracy scenarios $\left(\text{SCE}_{\sigma_y}\right)$.

In the following, \emph{quality of measurements} is used to qualitatively describe the difference between different measurement scenarios  $\left(\text{SCE}_{M}\right)$.
The \emph{worst quality} of measurements corresponds to data consisting of only the vapor mole fraction $y = [x^V_1]$ and perturbed using the larger standard deviations $\sigma^M_{x^V_1} = 0.001~\si{\mole\per\mole}$.
The \emph{best quality} corresponds to data consisting of both vapor mole fraction and temperature $y = [x^V_i, T]$ and perturbed using the smaller standard deviations $\sigma^M_{x^V_1}=0.0002~\si{\mole\per\mole},~\sigma^M_{T} = 0.01\si{K}$.
All scenarios are summarized in~\autoref{tab:meas-quality}.

\begin{table}
    \centering
    \caption{Summary of all measurement scenarios $\text{SCE}_{M}$. ``Worst'' and ``best'' refer to measurement quality. Non-diagonal scenarios are investigated, but have no label.}
    \label{tab:meas-quality}
    \begin{tabular}{cc|cc}
        \toprule
         & $\text{SCE}_{\sigma_y}$ & default & precise \\
        & $[\sigma^M_{x^V_1}~\si[per-mode=fraction]{\mole\per\mole},~\sigma^M_{T}~\si{K}]$ & $[0.001,~0.03]$ & $[0.0002,~0.01]$\\
       \midrule
        \multirow{2}{*}{$\text{SCE}_{y}$}& $y = [x^V_i]   $   & worst & - \\
    & $y = [x^V_i, T]$  & - & best \\
        \bottomrule
    \end{tabular}

\end{table}

Simulation and optimization are formulated and solved using the open-source Python package \emph{mopeds}\footnote{\href{https://pypi.org/project/mopeds/}{https://pypi.org/project/mopeds/}}~\citep{Kozachynskyi2024}.
It provides a top-level interface to a modeling library with automatic differentiation \emph{CasADi}~\citep{Andersson2019}.
Both the PE and OED optimization problems were formulated as nonlinear programming (NLP) programs and solved using IPOPT~\citep{Waechter2005} with MA57\footnote{HSL, a collection of Fortran codes for large-scale scientific computation. See \href{https://www.hsl.rl.ac.uk}{https://www.hsl.rl.ac.uk}}~\citep{Duff2004} as linear solver.

The parameters in the PE problem are bounded as follows: $-100 \leq A \leq 100$ and $-1.5\text{E}5 \leq B \leq 1.5\text{E}5$.
The parameter $\alpha$ is bounded between [0, 2] for the first case study and [0.1, 0.6] for the second case study.
As shown later, larger $\alpha$ bounds were not used for the second case study due to unstable performance of sOED-PE.
OED was solved using a multistart method, where 21 different combinations of $x^L_1$ and $P$ were uniformly distributed between optimizer bounds and used as a guess.

\subsection{Case Study I: FFD}

The results shown below are generated for each selected binary mixture using the $\left(\text{MC}^{\text{I}}\right)$ in~\autoref{MCi}:

\begin{enumerate}
    \item Number of Monte Carlo simulations $N_{MC}=1000$
    \item The measurement grid $U$ consists of 20 liquid mole fractions $x^L_1$ uniformly distributed between 0.01 and 0.99~\si{\mole\per\mole} at pressures $P$ 0.5 and 1.5~\si{bara}.
    \item The prediction grid $U^p$ consists of 20 liquid mole fractions $x^L_1$ uniformly distributed between 0.01 and 0.99~\si{\mole\per\mole} at pressure $P$ 1~\si{bara}.
    \item 16 different parameter scenarios $\left(\text{SCE}_{\theta}\right)$ are tested for each mixture:
        \begin{itemize}
            \item \emph{All}: all five NRTL parameters are considered identifiable, no regularization is used.
            \item $A^\ast_{12},B^\ast_{12},\ldots,\alpha^\ast$: one of the five parameters is fixed to its correct value and is not estimated by PE.
            \item $\underline{A_{12}},\overline{A_{12}},\ldots,\overline{B_{21}}$: one of the four parameters is fixed to its wrong value generated by perturbation of the correct value, e.g.,  $\underline{A_{12}} = 0.8 \cdot A^\ast_{12}$ or $\overline{A_{12}} = 1.2 \cdot A^\ast_{12}$.
            \item $\underline{\alpha},\overline{\alpha}$: the parameter is fixed to the wrong value $\underline{\alpha}=0.1$ or $\overline{\alpha}=0.6$.
        \end{itemize}
    \item 4 different measurement scenarios $\text{SCE}_{M}$ are tested for each mixture and parameter scenario $\text{SCE}_{\theta}$ combined from two measurement accuracy scenarios $\text{SCE}_{\sigma_y}$ and two measurement variable scenarios $\text{SCE}_{y}$ shown in~\autoref{tab:meas-quality}.
\end{enumerate}

In total, 732 different scenarios were analyzed, summarized in~\autoref{tab:cs-scenarios}.

\begin{table}
\centering
\caption{Overview of selected scenarios (SCE) for different mixtures used in \emph{Case Study I} combined from measurement scenarios $\text{SCE}_M$ and parameter scenarios $\text{SCE}_{\theta}$. Special cases: $^a~A_{12}=A_{21}=0$; $^b~\underline{\alpha} = \alpha^\ast = 0.1$.}
\label{tab:cs-scenarios}
\begin{tabular}{lcccc}
\toprule
 & $\text{SCE}_{M}$ & $\text{SCE}_{\theta}$ & Total \\
Mixture &  &  &  \\
\midrule
Acetone / Chloroform & 4 & 16 & 64 \\
Methanol / Benzene & 4 & 16 & 64 \\
Methanol / Water & 4 & 16 & 64 \\
Benzene / Hexafluorobenzene & 4 & $12^a$ & 48 \\
Chloroform / Ethyl acetate & 4 & $12^a$ & 48 \\
Ethanol / Benzene & 4 & 16 & 64 \\
Ethanol / Dioxane & 4 & 16 & 64 \\
Acetone / Ethanol & 4 & 16 & 64 \\
Methanol / Ethanol & 4 & 16 & 64 \\
Methanol / Epoxybutane & 4 & $15^b$ & 60 \\
Chloroform / Epoxybutane & 4 & 16 & 64 \\
Ethanol / Acetic acid & 4 & 16 & 64 \\
\midrule
Total &  &  & 732 \\
\bottomrule
\end{tabular}
\end{table}

\subsubsection{Assessment of linearity assumption}

In general, we observed no correlation between the type of azeotropic behavior of the mixture and the parameter linearity of the NRTL model.
Even though it might be possible to find some correlation between parameter uncertainty, nonlinearity, and mixture properties, e.g., ratio of the vapor pressures, or azeotropic behaviour, it has become evident that such analysis might require more binary mixtures.

\begin{figure*}
    \centering
    \includegraphics[width=200pt]{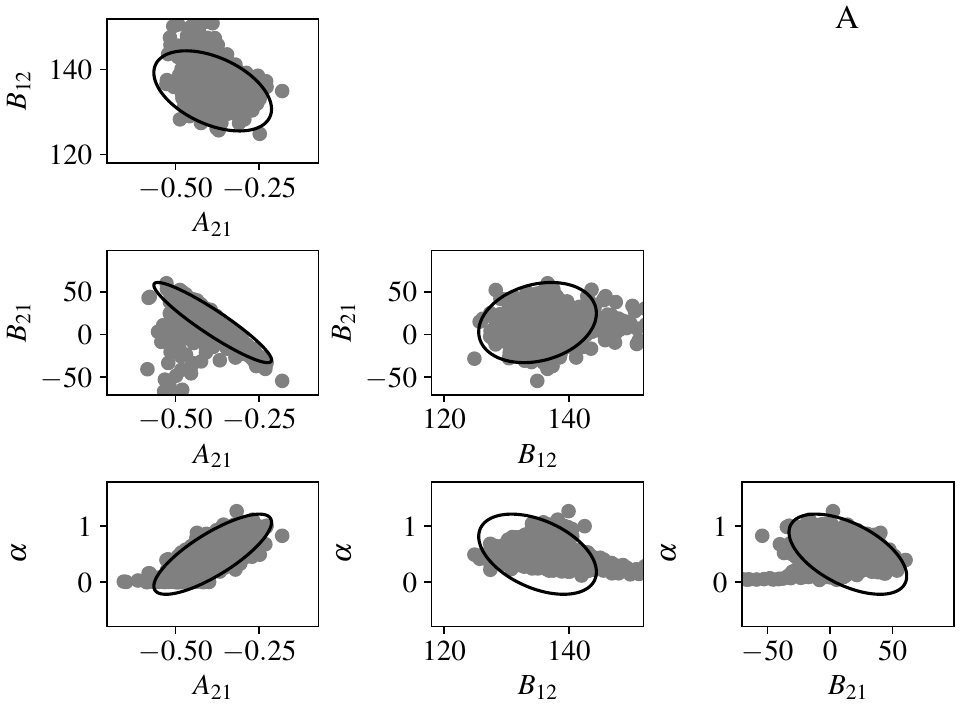}
    \includegraphics[width=200pt]{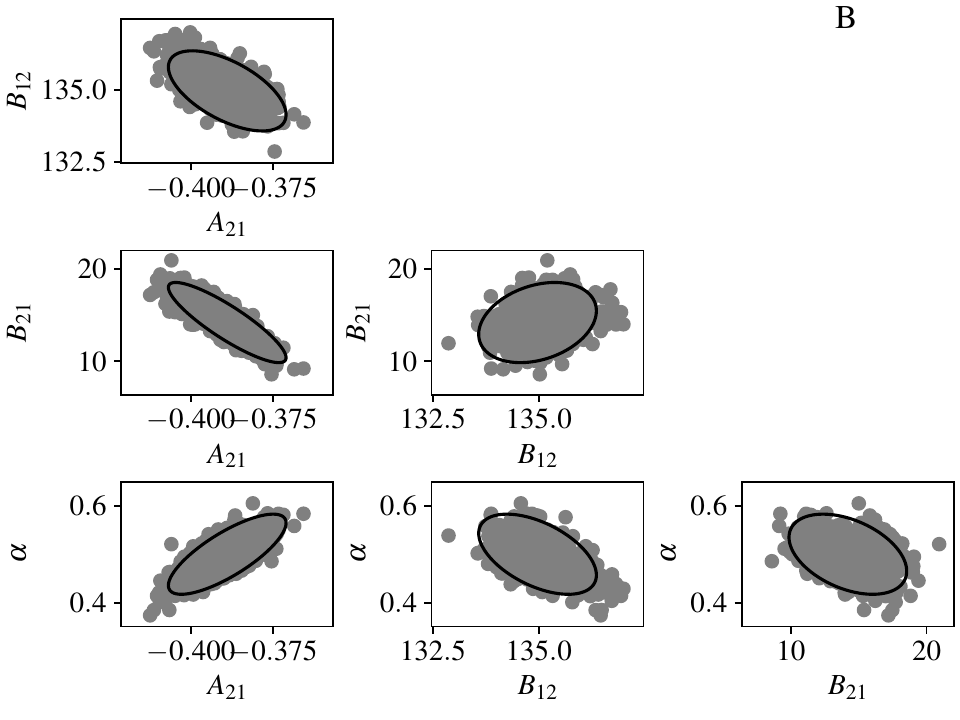}
    \caption{Parameter covariance plots for Ethanol / Acetic Acid mixture at worst (A) and best (B) measurement quality $\text{SCE}_{M}$ for the scenario $\text{SCE}_{\theta}=A^\ast_{12}$: gray markers represent values computed using Monte Carlo simulations, the line represents the 95\% confidence ellipse of the true parameter $\theta^\ast$ computed using linearization assumption.}
    \label{fig:covariance-B-fixed}
\end{figure*}

As expected, the parameter nonlinearity decreases with increasing measurement quality.
In~\autoref{fig:covariance-B-fixed} the parameter uncertainty calculated with the parameter $A_{12}$ fixed to its ``true'' value $\left(\text{SCE}_{\theta} = A_{12}^{\ast}\right)$ is visualized for the best and worst measurement quality scenarios $\text{SCE}_{M}$~(\autoref{tab:meas-quality}).
At the best measurement quality, the model is less nonlinear in parameters, which can be clearly seen as the parameter uncertainty follows the shape of the confidence ellipse calculated using the linearization assumption.
On the other hand, at the worst measurement quality, the parameters calculated using Monte Carlo simulations do not form the ellipse, so the model is more nonlinear in parameters.

The average values of the metrics for all mixtures shown in~\autoref{tab:cs-stats}, are used for further analysis.
The following trend was observed for all mixtures: fixing the parameter $\alpha$ to the correct value $\left(\text{SCE}_{\theta} = \alpha^{\ast}\right)$ always made the covariance matrix more elliptical and the effect of parameters on the model less nonlinear, as indicated by high values of the parameter normality $\overline{W}$ and the model linearity $Q^{0.95}_{\theta}$.
This consistent trend cannot be observed for the other parameters $A$ and $B$.

\begin{table}
\centering
\caption{Average and standard deviation of parameter normality $\overline{W}$, linearity $Q^{0.95}_{\theta}$, and prediction quality $Q^{0.95}_y$ across all mixtures and measurements scenarios $\text{SCE}_{M}$ for different parameter scenarios $\text{SCE}_{\theta}$.}
\label{tab:cs-stats}
\begin{tabular}{cccccc}
\toprule
                      & $\overline{W}$  & $Q^{0.95}_{\theta}$ & $Q^{0.95}_y$     \\
$\text{SCE}_{\theta}$ &                 &                     &                  \\
\midrule                                                                        
\emph{All}            & 0.86 $\pm$ 0.22 & 0.97 $\pm$ 0.02     & 0.95 $\pm$ 0.01  \\
\midrule                                                                        
$A^\ast_{12}$         & 0.93 $\pm$ 0.19 & 0.97 $\pm$ 0.02     & 0.95 $\pm$ 0.02  \\
$A^\ast_{21}$         & 0.93 $\pm$ 0.22 & 0.98 $\pm$ 0.01     & 0.95 $\pm$ 0.02  \\
$B^\ast_{12}$         & 0.93 $\pm$ 0.20 & 0.97 $\pm$ 0.02     & 0.95 $\pm$ 0.01  \\
$B^\ast_{21}$         & 0.94 $\pm$ 0.16 & 0.97 $\pm$ 0.03     & 0.95 $\pm$ 0.04  \\
$\alpha^\ast$         & 1.00 $\pm$ 0.00 & 0.98 $\pm$ 0.00     & 0.96 $\pm$ 0.00  \\
\midrule                                                                        
$\overline{A_{12}}$   & 0.91 $\pm$ 0.21 & 0.65 $\pm$ 0.16     & 0.67 $\pm$ 0.29  \\
$\underline{A_{12}}$  & 0.84 $\pm$ 0.29 & 0.46 $\pm$ 0.19     & 0.62 $\pm$ 0.28  \\
$\overline{A_{21}}$   & 0.86 $\pm$ 0.26 & 0.61 $\pm$ 0.18     & 0.66 $\pm$ 0.29  \\
$\underline{A_{21}}$  & 0.92 $\pm$ 0.20 & 0.56 $\pm$ 0.14     & 0.66 $\pm$ 0.27  \\
$\overline{B_{12}}$   & 0.86 $\pm$ 0.27 & 0.62 $\pm$ 0.15     & 0.71 $\pm$ 0.28  \\
$\underline{B_{12}}$  & 0.91 $\pm$ 0.19 & 0.57 $\pm$ 0.17     & 0.67 $\pm$ 0.26  \\
$\overline{B_{21}}$   & 0.88 $\pm$ 0.23 & 0.72 $\pm$ 0.12     & 0.59 $\pm$ 0.30  \\
$\underline{B_{21}}$  & 0.83 $\pm$ 0.30 & 0.56 $\pm$ 0.23     & 0.51 $\pm$ 0.34  \\
$\overline{\alpha}$   & 0.94 $\pm$ 0.17 & 0.53 $\pm$ 0.16     & 0.24 $\pm$ 0.24  \\
$\underline{\alpha}$  & 1.00 $\pm$ 0.00 & 0.57 $\pm$ 0.15     & 0.34 $\pm$ 0.31  \\
\bottomrule
\end{tabular}
\end{table}

\begin{figure}[h]
    \centering
    \includegraphics[width=\columnwidth]{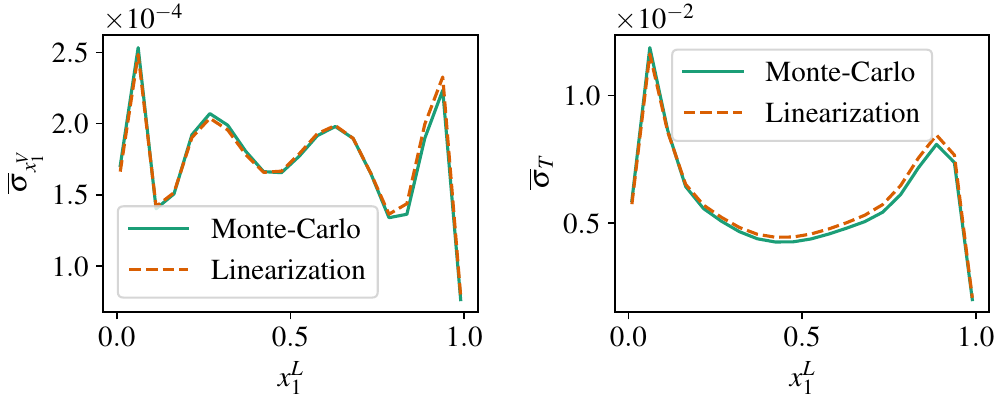}
    \caption{Standard deviation of the model prediction for Ethanol / Benzene mixture: real and estimated using linearization assumption. Default precision scenario $\text{SCE}_{\sigma_y}$, all parameters are estimated $\text{SCE}_\theta = All$.} 
    \label{fig:goodpredictionaccuracy}
\end{figure}

The model prediction quality $Q^{0.95}_y$ is around the expected value of 0.95 for the parameter scenarios where either all parameters were estimated $\left(\text{SCE}_{\theta} = All\right)$, or one of the parameters was fixed to the ``correct'' value $\left(\text{SCE}_{\theta} = \theta^\ast\right)$~(\autoref{tab:cs-stats}).
It shows that the prediction accuracy of the model approximated with linearization assumptions~(\autoref{eq:prediction-std}) is correct, or at least a valid underestimation of the real prediction accuracy, for all analyzed mixtures, despite parameters showing nonlinear effects, as indicated by the average and standard deviation values of the parameter normality $\overline{W}$.
The similarity of the real and linearly approximated standard deviations of the model prediction can be clearly observed in~\autoref{fig:goodpredictionaccuracy}.

Usually, when the parameter uncertainty of the estimated parameters is published, if at all, authors only show confidence intervals of the individual parameters.
For example, the individual confidence intervals of the parameters used to estimate the prediction accuracy in~\autoref{fig:goodpredictionaccuracy} are shown in~\autoref{tab:covariance}.
Since the uncertainty of the parameter $B_{12}$ is in the same order of magnitude as the parameter value itself (\autoref{tab:covariance}), the parameter estimate looks unreliable.
So how can the prediction accuracy shown in~\autoref{fig:goodpredictionaccuracy} be accurate even though the corresponding parameters look uncertain at first glance~(\autoref{tab:covariance})?
It depends on the fact that colinearity, the dependence of the parameters on each other, plays a crucial role.
Colinearity is described by the non-diagonal elements of the covariance matrix shown in~\autoref{tab:parcovmatrix}.
Since the individual parameter confidence intervals do not include the effect of colinearity, we recommend publishing only the parameter variance--covariance matrix~(\autoref{tab:covariance}) with the corresponding prediction accuracy~(\autoref{fig:goodpredictionaccuracy}) for the estimated NRTL parameters.

\begin{table}
    \centering
    \caption{Estimated parameter values with individual confidence intervals. Default precision scenario $\text{SCE}_{\sigma_y}$, all parameters are estimated $\text{SCE}_\theta = All$.}
\label{tab:covariance}
    \begin{tabular}{cccc}
    \toprule
    Parameter & Value & CI 95\%, absolute & CI 95\%, relative \\
    \midrule
    $A_{12}$ & 0.568 & $\pm$ 0.105 & $\quad$ $\pm$ 18.6~\%\\
    $A_{21}$ & -0.915& $\pm$ 0.151 & $\quad$ $\pm$ 16.5~\%\\
    $B_{12}$ & -54.8 & $\pm$ 40.93 & $\quad$ $\pm$ 74.6~\%\\
    $B_{21}$ & 882.0 & $\pm$ 52.50 & $\quad$ $\pm$ 5.90~\%\\
    $\alpha$ & 0.3   & $\pm$ 0.014 & $\quad$ $\pm$ 4.80~\%\\
    \bottomrule
    \end{tabular}
\end{table}

\begin{table}
    \centering
    \caption{Parameter variance--covariance matrix. Default precision scenario $\text{SCE}_{\sigma_y}$, all parameters are estimated $\text{SCE}_\theta = All$.}
    \label{tab:parcovmatrix}
\begin{tabular}{l|ccccc}
\toprule
 & $A_{12}$ & $A_{21}$ & $B_{12}$ & $B_{21}$ & $\alpha$ \\
\midrule
$A_{12}$ & 2.7\texttimes$10^{-3}$ & -3.5\texttimes$10^{-3}$ & -1.01 & 1.27 & -1.1\texttimes$10^{-4}$ \\
$A_{21}$ & -3.5\texttimes$10^{-3}$ & 0.01 & 1.24 & -1.93 & 3.3\texttimes$10^{-5}$ \\
$B_{12}$ & -1.01 & 1.24 & 409 & -459 & 0.08 \\
$B_{21}$ & 1.27 & -1.93 & -459 & 672 & -0.04 \\
$\alpha$ & -1.1\texttimes$10^{-4}$ & 3.3\texttimes$10^{-5}$ & 0.08 & -0.04 & 5.2\texttimes$10^{-5}$ \\
\bottomrule
\end{tabular}
\end{table}

\subsubsection{Consequences of fixing $\alpha$}

\begin{figure*}
    \centering
    \includegraphics[width=200pt]{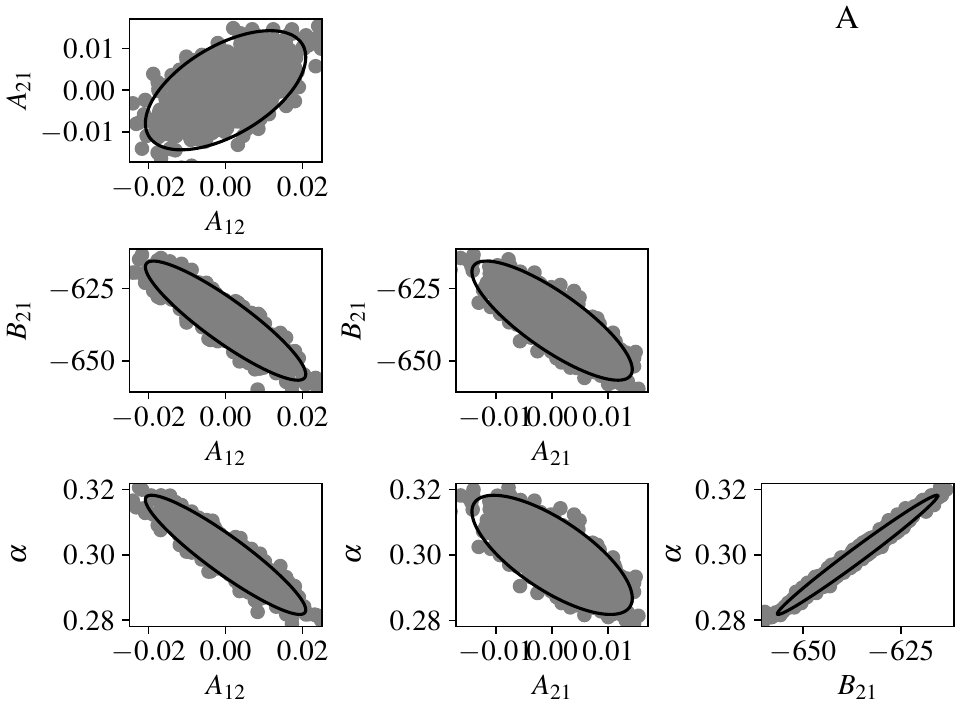}
    \includegraphics[width=200pt]{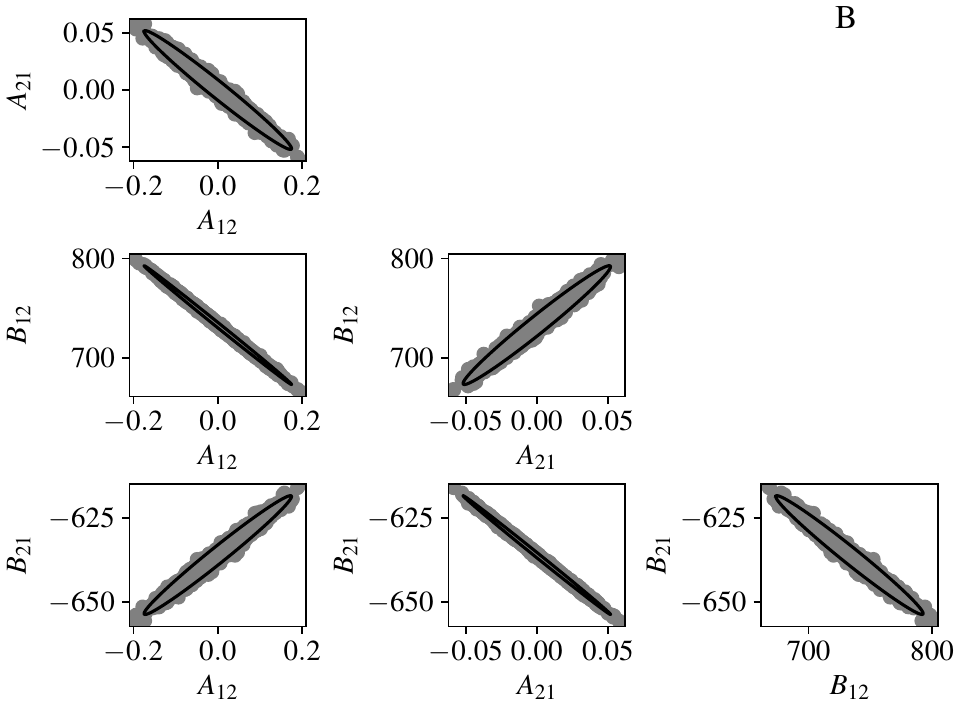}
    \caption{Parameter covariance plots for Chloroform / Ethyl Acetate mixture with different fixed parameters: $B_{12}^{\ast}$ (A) and $\alpha^{\ast}$ (B) for best measurement quality scenario $\text{SCE}_{M}$. Gray markers represent values computed using the Monte Carlo method, the line represents the 95\% confidence ellipse of the true parameter $\theta^\ast$ computed using linearization assumption.}
    \label{fig:covariance-chlor-epoxy}
\end{figure*}

As mentioned earlier, fixing $\alpha$ consistently leads to a model that is more linear in parameters.
However, an important drawback of fixing $\alpha$ has been noticed: it can lead to identifiaibility issues.
This can be seen in~\autoref{fig:covariance-chlor-epoxy}, where the parameter covariance is generated for two parameter scenarios: either $B_{12}$ $\left(\text{SCE}_{\theta} = B_{12}^{\ast}\right)$ or $\alpha$ $\left(\text{SCE}_{\theta} = \alpha^{\ast}\right)$ were fixed to their correct values.
In both cases, the model is less nonlinear in parameters, i.e., the exact pairwise parameter uncertainty has the shape of an ellipse.
However, the identifiability problem is clearly seen for the $\alpha^\ast$ scenario: the ellipses are much more elongated, indicating that the remaining (non-fixed) parameters are highly correlated.

So far we have only analyzed the scenarios where $\alpha$ is fixed to its ``true'' value.
The consequence of fixing $\alpha$ to incorrect values is clearly shown in~\autoref{tab:cs-stats}: for ``wrong`` scenarios $\left(\text{SCE}_{\theta} = \left[\overline\alpha, \underline\alpha\right]\right)$ the prediction quality $Q^{0.95}_{y}$ is very low, which means that PE with wrongly fixed $\alpha$ leads to biased models with bad prediction quality.

The question, of course, is how to determine in practice whether the model is biased.
In~\autoref{tab:table-fixed-params-comparison} the metrics for a Methanol / Water mixture are shown for different parameter scenarios $\text{SCE}_{\theta}$.
Based on the low value of the prediction quality metric $Q^{0.95}_y$ it can be concluded that the model is biased for both scenarios $\text{SCE}_{\theta} = \left[\underline\alpha, \overline\alpha\right]$.
However, since the prediction quality metric cannot be calculated without knowing the ``true'' parameter values, one would either use the validation data set or check the estimation accuracy $s_{x_1^V}$. 
As described earlier, it should be close to the measurement accuracy $\sigma^M_{x^V_i}$.
While the estimation accuracy in both $\text{SCE}_{\theta} = \overline\alpha$ and $\text{SCE}_{\theta} = \underline\alpha$ scenarios is in the same order of magnitude as the measurement accuracy, the experimenter might conclude that the model is of good quality.
This highlights the need for caution when manually fixing parameters, as parameters that appear good at first glance may lead to incorrect model predictions.

\begin{table}
\centering
\caption{Prediction quality $Q^{0.95}_y$, estimation accuracy $\overline{s}_{x_1^V}$, parameter normality $\overline{W}$, prediction accuracy $\overline{\sigma}_{x_1^V}$, and parameter linearity $Q^{0.95}_{\theta}$ for the Methanol / Water mixture at best measurement quality $\text{SCE}_{M}$, i.e., $\sigma^M_{x^V_1}=0.0002$, for different parameter scenarios $\text{SCE}_{\theta}$.}
\label{tab:table-fixed-params-comparison}
\begin{tabular}{cccccc}
\toprule
                        & $Q^{0.95}_y$ & $\overline{s}_{x_1^V}$ & $\overline{W}$ & $\overline{\sigma}_{x_1^V}$ & $Q^{0.95}_{\theta}$ \\
$\text{SCE}_{\theta}$   &              &                        &                &                             &                     \\
\midrule                                                        
\emph{All}              & 0.96         & 1.9\texttimes10$^{-4}$ & 0.92           & 3.3\texttimes10$^{-5}$      & 0.97                \\
\midrule                                                        
$A^\ast_{12}$           & 0.96         & 1.9\texttimes10$^{-4}$ & 1.00           & 2.6\texttimes10$^{-5}$      & 0.98                \\
$A^\ast_{21}$           & 0.96         & 1.9\texttimes10$^{-4}$ & 1.00           & 2.6\texttimes10$^{-5}$      & 0.98                \\
$B^\ast_{12}$           & 0.95         & 1.9\texttimes10$^{-4}$ & 0.99           & 2.6\texttimes10$^{-5}$      & 0.97                \\
$B^\ast_{21}$           & 0.96         & 1.9\texttimes10$^{-4}$ & 1.00           & 2.6\texttimes10$^{-5}$      & 0.98                \\
$\alpha^\ast$           & 0.96         & 2.0\texttimes10$^{-4}$ & 1.00           & 2.7\texttimes10$^{-5}$      & 0.98                \\
\midrule                                                        
$\overline{A_{12}}$     & 0.71         & 3.2\texttimes10$^{-4}$ & 0.24           & 2.7\texttimes10$^{-5}$      & 0.74                \\
$\underline{A_{12}}$    & 0.25         & 4.0\texttimes10$^{-4}$ & 0.11           & 2.6\texttimes10$^{-5}$      & 0.25                \\
$\overline{A_{21}}$     & 0.92         & 1.9\texttimes10$^{-4}$ & 1.00           & 2.6\texttimes10$^{-5}$      & 0.50                \\
$\underline{A_{21}}$    & 0.90         & 2.0\texttimes10$^{-4}$ & 1.00           & 2.6\texttimes10$^{-5}$      & 0.50                \\
$\overline{B_{12}}$     & 0.70         & 2.0\texttimes10$^{-4}$ & 1.00           & 2.6\texttimes10$^{-5}$      & 0.50                \\
$\underline{B_{12}}$    & 0.28         & 2.2\texttimes10$^{-4}$ & 0.99           & 2.6\texttimes10$^{-5}$      & 0.50                \\
$\overline{B_{22}}$     & 0.91         & 2.0\texttimes10$^{-4}$ & 1.00           & 2.6\texttimes10$^{-5}$      & 0.75                \\
$\underline{B_{22}}$    & 0.87         & 2.0\texttimes10$^{-4}$ & 1.00           & 2.6\texttimes10$^{-5}$      & 0.27                \\
$\overline{\alpha}$     & 0.14         & 2.8\texttimes10$^{-4}$ & 1.00           & 2.7\texttimes10$^{-5}$      & 0.50                \\
$\underline{\alpha}$    & 0.36         & 2.1\texttimes$10^{-4}$ & 1.00           & 2.7\texttimes10$^{-5}$      & 0.50                \\
\bottomrule
\end{tabular}
\end{table}

In summary, fixing $\alpha$ has long been a good solution to the obvious problem: estimating four parameters is easier than estimating five. 
However, the solution comes at the cost of a model potentially being biased because $\alpha$ is fixed to a wrong value.
The only alternative is to estimate all 5 NRTL parameters simultaneously, which has become an easier task due to advances in optimization solvers that perform their own regularization of the nonlinear programming problem, e.g., IPOPT.

\subsubsection{Regularization}

The previous analysis showed that manually fixing parameters might be a bad idea and that all 5 NRTL parameters should be estimated simultaneously.
Estimating all five parameters could lead to an ill-conditioned optimization problem that must be solved using regularization.
Such regularization - whether done manually, by subset selection, or automatically by an NLP solver - will introduce some bias that will manifest itself in a reduced prediction quality, which must be kept small.
But is it possible to use regularization without introducing significant bias?

For further analysis data was generated for $\text{SCE}_{\theta} = All$ scenario using method $\text{MC}^{\text{Reg}}$ in~\autoref{MCreg}.
The four regularization algorithms introduced before were applied $N_{MC}=1000$ times for all mixtures and measurement scenarios shown in~\autoref{tab:cs-scenarios}, resulting in 16 runs for every mixture.

In~\autoref{tab:cs-regularization-different} the average prediction accuracy $\overline{\sigma}_{x_1^V}$ and quality $Q^{0.95}_y$, as well as average frequency of the parameter marked as identifiable over all Monte Carlo simulations are shown for one mixture and one measurement scenario $\text{SCE}_{M}$.
It can be seen that after estimation of all five parameters, the \emph{GO} regularization algorithm selected both parameters $B_{21}, A_{21}$ as identifiable in only $\approx65\%$ of the time.
In the remaining 35\%, these parameters were fixed, resulting in an average prediction accuracy improvement of 15\% while keeping the prediction quality $Q^{0.95}_y$ high.
This means that while the \emph{GO} algorithm still introduces bias, as any regularization does, the resulting model is not significantly biased on average.
In contrast, the more conservative \emph{FS} algorithm suggests fixing more parameters on average, resulting in a more biased model $Q^{0.95}_y\approx0.9$ on average.

\begin{table}
\centering
\caption{Prediction accuracy $\overline{\sigma}_{x_1^V}$, prediction quality $Q^{0.95}_y$ and identifiability of parameters for different regularization scenario $\text{SCE}_{Reg}$. Rows with parameter names indicate the frequency with which the parameter is identifiable, where 1 means that the parameter was identifiable in all Monte Carlo simulations. Results are shown for the Acetone / Chloroform mixture, with data generated using the worst measurement quality $\text{SCE}_{M}$.}
\label{tab:cs-regularization-different}
\begin{tabular}{c|ccccc}
\toprule
& \multicolumn{4}{c}{$\text{SCE}_{Reg}$} \\
& E & GO & SVD & FS \\
\midrule
$\overline{\sigma}_{x_1^V}$ & 2.7\texttimes$10^{-4}$ & 2.3\texttimes$10^{-4}$ & 2.2\texttimes$10^{-4}$ & 2.0\texttimes$10^{-4}$ \\
\noalign{\smallskip}
$Q^{0.95}_y$ & 0.97 & 0.96 & 0.94 & 0.90 \\
\midrule
$A_{12}$ & 0.93 & 0.88 & 0.87 & 0.82 \\
$A_{21}$ & 0.87 & 0.69 & 0.68 & 0.43 \\
$B_{12}$ & 1.00 & 0.98 & 0.98 & 0.98 \\
$B_{21}$ & 0.83 & 0.61 & 0.58 & 0.31 \\
$\alpha$ & 1.00 & 0.94 & 0.95 & 0.50 \\
\bottomrule
\end{tabular}
\end{table}

The overall performance of the four regularization techniques is summarized in~\autoref{tab:reg-stats-overall}, while the influence of different measurement scenarios is presented in~\autoref{tab:reg-statistics} and the frequency of identifiable parameters in~\autoref{tab:ident-parameter}.
The regularization techniques can be ranked based on the average number of identifiable parameters in descending order as follows: \emph{E}, \emph{GO}, \emph{SVD}, and \emph{FS}.
Generalized orthogonalization \emph{GO} proved itself as the most effective algorithm: it regularized most of the time without affecting the prediction accuracy and without any thresholds.
Algorithms \emph{E} and \emph{FS} have not been tuned and may perform better after tuning, but there is no obvious methodology for doing so.
The \emph{SVD} algorithm with a condition number threshold of 1000 was too conservative with respect to the number of identifiable parameters, resulting in a bias in some cases.

\begin{table}
    \centering
    \caption{Average and standard deviation of prediction quality $Q^{0.95}_y$ and number of identifiable parameters $\overline{N}_i$ for 4 regularization techniques $\text{SCE}_{Reg}$.}
    \label{tab:reg-stats-overall}
\begin{tabular}{ccc}
\toprule
 & $Q^{0.95}_y$ & $\overline{N}_i$ \\
$\text{SCE}_{Reg}$ &  &  \\
\midrule
E & 0.95 $\pm$ 0.04 & 4.76 $\pm$ 0.59 \\
GO & 0.94 $\pm$ 0.04 & 4.43 $\pm$ 0.67 \\
SVD & 0.91 $\pm$ 0.04 & 3.93 $\pm$ 0.78 \\
FS & 0.87 $\pm$ 0.04 & 2.96 $\pm$ 0.82 \\
\bottomrule
\end{tabular}
\end{table}

\subsection{Case Study II: sOED-PE}

In the second case study, the results are generated using the methods $\text{MC}_{sOED-PE}^{\text{I}}$ and $\text{MC}_{sOED-PE}^{\text{Reg}}$ from ~\autoref{MCsoedi} and~\autoref{MCsoedreg} with the following settings:

\begin{enumerate}
    \item Number of repeated sOED-PE runs $N_{MC} = 500$
    \item Number of sOED-PE iterations $N_{sOED-PE}=15$
    \item Initial measurement grid consists of the following points $U^{init}=[x_{i=1}, P$]: [0.05, 0.5e5], [0.95, 0.5e5], [0.05, 1.5e5], [0.95, 1.5e5], [0.5, 1e5], [0.65, 1e5].
    \item Both measurement variables are used $\text{SCE}_{y}=[x^V_i, T]$.
    \item The A criterion is used for OED.
    \item 6 different regularization scenarios $\text{SCE}_{Reg}$ are tested:
        \begin{enumerate}
            \item No regularization is used.
            \item 4 regularization techniques are used for PE only.
            \item \emph{GO} technique is used for both PE and OED.
        \end{enumerate}
    \item Three different parameter scenarios $\text{SCE}_{\theta}$ are tested: $\alpha^\ast,~\overline\alpha$, and $\underline\alpha$.
    \item All mentioned scenarios are tested for each mixture in~\autoref{tab:cs-scenarios} and two different measurement accuracies $\text{SCE}_{\sigma_y}$(~\autoref{tab:meas-quality}).

\end{enumerate}

A total of 214 different scenarios were executed, 18 for every mixture, as shown in~\autoref{tab:cs2-scenarios}.
In general, the sOED-PE performance is relatively robust for most scenarios, typically less than 10\% of all Monte Carlo runs failed at some point during the sOED-PE iterations.
Data from failed runs are excluded from further analysis.

The $\alpha$ bounds for PE had to be set between 0.1 and 0.6.
Although it would be interesting to analyze the results of sOED-PE for larger $\alpha$ bounds, this was not possible due to the high number of failing runs.
Failures were usually caused by non-convergence / local infeasibility of either PE or OED optimization problems.

\subsubsection{sOED-PE performance}

The average metrics for all mixtures are shown for different parameter scenarios $\text{SCE}_{\theta}$ and the default measurement accuracy scenario $\text{SCE}_{\sigma_y}$ in~\autoref{tab:soedpe-stats}.
Overall, as expected, the prediction accuracy $\overline{\sigma}_{x^V_1}$ improves as the number of sOED-PE interations increases.
The estimation accuracy $\overline{s}_{x^V_1}$ for all scenarios is close to the measurement error $\sigma^M_{x^V_i}=0.001$, which is not expected for ``wrong'' scenarios $\text{SCE}_{\theta}=\left[\underline\alpha, \overline\alpha\right]$.
This will be analyzed in detail in the next section.

The prediction quality $Q^{0.95}_y$ worsens over the sOED-PE iterations and becomes less than 0.95, indicating that increased prediction accuracy comes at the expense of model prediction bias.
This is due to the fact that only the parameter uncertainty is included in the A criterion for OED and not the prediction accuracy.
This means that OED only suggests experiments in a few regions of the design space, as will be shown later,  and as the number of experiments increases, the model becomes ``overfitted'' to the regions of the selected design space.
This behavior is expected for all the classical OED criteria based only on parameter accuracy, i.e., A, D, E, etc.

\begin{table}
    \centering
    \caption{Average prediction accuracy $\overline{\sigma}_{x_1^V}$ and estimation accuracy  $\overline{s}_{x_1^V}$, and average and standard deviation of prediction quality $Q^{0.95}_{\sigma_{x^V_1}}$ for all mixtures for the first and last sOED-PE iteration $I_{sOED-PE}$. Discrepancy $D^{U}$ is calculated for all experiments. Different parameter scenarios $\text{SCE}_{\theta}$. Default measurement accuracy scenario $\text{SCE}_{\sigma_y}$.}
    \label{tab:soedpe-stats}
\begin{tabular}{cc|cccc}
\toprule
\rot{$\text{SCE}_{\theta}$}           & \rot{$I_{sOED-PE}$} & $\overline{\sigma}_{x_1^V}$ & $\overline{s}_{x_1^V}$ & $Q^{0.95}_{x_1^V}$ & $D^{U}$ \\
\midrule                                                                                                           
\multirow{3}{*}{$\alpha^\ast$}        & 1                   & 3.8\texttimes$10^{-4}$      & 1.0\texttimes$10^{-3}$ & 0.99 $\pm$ 0.05    & - \\
                                      & 15                  & 1.5\texttimes$10^{-4}$      & 9.9\texttimes$10^{-4}$ & 0.96 $\pm$ 0.11    & - \\
                                      & All                 & -                           & -                      & -                  & 6.7\texttimes$10^{-2}$ \\
\cline{1-6}                                                                                                        
\multirow{3}{*}{$\underline{\alpha}$} & 1                   & 4.0\texttimes$10^{-4}$      & 1.0\texttimes$10^{-3}$ & 0.94 $\pm$ 0.15    & - \\
                                      & 15                  & 1.4\texttimes$10^{-4}$      & 1.0\texttimes$10^{-3}$ & 0.74 $\pm$ 0.28    & - \\
                                      & All                 & -                           & -                      & -                  & 4.9\texttimes$10^{-2}$ \\
\cline{1-6}                                                                                                        
\multirow{3}{*}{$\overline{\alpha}$}  & 1                   & 4.0\texttimes$10^{-4}$      & 1.0\texttimes$10^{-3}$ & 0.93 $\pm$ 0.13    & - \\
                                      & 15                  & 1.5\texttimes$10^{-4}$      & 1.0\texttimes$10^{-3}$ & 0.72 $\pm$ 0.25    & - \\
                                      & All                 & -                           & -                      & -                  & 6.6\texttimes$10^{-2}$ \\
\cline{1-6}                                                                                                        
\multirow{3}{*}{\emph{All}}           & 1                   & 5.5\texttimes$10^{-4}$      & 1.0\texttimes$10^{-3}$ & 1.00 $\pm$ 0.02    & - \\
                                      & 15                  & 2.0\texttimes$10^{-4}$      & 1.0\texttimes$10^{-3}$ & 0.97 $\pm$ 0.09    & - \\
                                      & All                 & -                           & -                      & -                  & 3.1\texttimes$10^{-2}$ \\
\cline{1-6}
\bottomrule
\end{tabular}
\end{table}

\begin{figure}[h]
    \centering
    \includegraphics[width=\columnwidth]{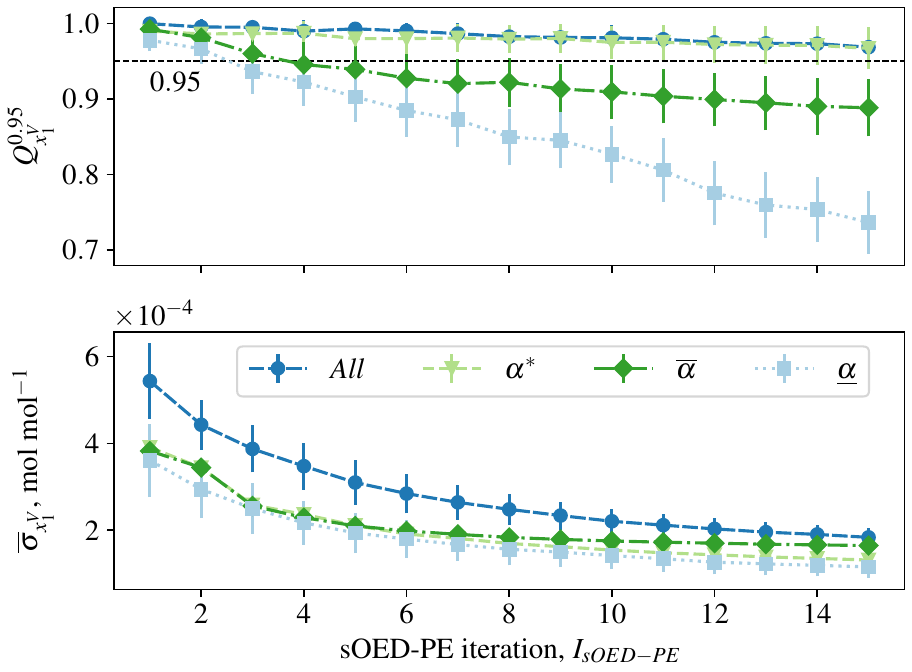}
    \caption{Change of average prediction quality $Q^{0.95}_y$ and accuracy $\overline{\sigma}_{x_i^V}$ over sOED-PE iterations $I_{sOED-PE}$ for different scenarios. Acetone / Chloroform mixture, both measurement variables are used $\text{SCE}_{y}$ with default measurement accuracy $\text{SCE}_{\sigma_y}$. Standard deviation of the prediction quality $Q^{0.95}_y$ is scaled down by a factor of 0.25 for better visualization.}
    \label{fig:wrong-alpha-soed}
\end{figure}

The change in prediction quality and prediction accuracy is shown in~\autoref{fig:wrong-alpha-soed} for one mixture at the default measurement accuracy scenario $\text{SCE}_{\sigma_y}$~(\autoref{tab:meas-quality}).
The importance of using the sOED-PE method is clearly shown for the scenario $\text{SCE}_{\theta}=All$, as the prediction accuracy decreases with each new experiment performed, so the model becomes more accurate.
Compared to the standard FFD approach, where all experiments are performed at once, the use of sOED-PE methods allows to constantly monitor the prediction accuracy and to make a decision: either to run a new experiment or to stop the experimental campaign because the target accuracy has already been reached.

\subsubsection{Consequences of fixing $\alpha$}

\begin{figure}[h]
    \centering
    \includegraphics[width=\columnwidth]{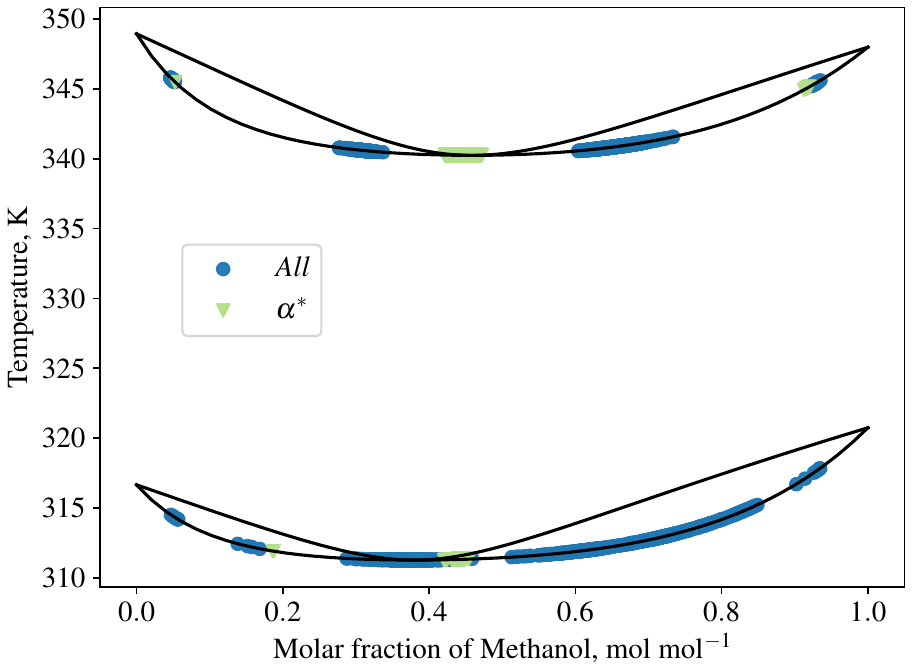}
    \caption{Experimental design of sOED-PE for different scenarios. Black lines are the predicted vapor and liquid molar fractions at 0.5 and 1.5~\si{bara}. Methanol / Epoxybutane mixture, best measurement quality scenario $\text{SCE}_{M}$.}
    \label{fig:soed-grid}
\end{figure}

As in the previous case study, fixing $\alpha$ to ``wrong'' values mostly leads to a biased model, as shown in~\autoref{tab:soedpe-stats}.
The model bias is difficult to detect by the estimation accuracy metric $\overline{s}_y$ when using sOED-PE.
The reason for this can be clearly seen in the difference in the final experimental design shown in~\autoref{fig:soed-grid} when either all 5 $\left(\text{SCE}_{\theta}=All\right)$ or only 4 $\left(\text{SCE}_{\theta}=\alpha^\ast\right)$ NRTL parameters are used.
With a smaller number of parameters ($\text{SCE}=\alpha^\ast$), OED with the A criterion as objective function suggests executing experiments mostly at 3 distinct points: near azeotropes at the lowest and the highest pressure levels, and at higher $x^L_1$ concentration at the highest pressure.
Of course, in this case the estimation accuracy $\overline{s}_y$ will be close to the measurement uncertainty $\sigma_y^M$ because technically four parameters are estimated using only three different experimental conditions.
The same effect can be seen for all mixtures as the discrepancy metric $D^{U}$ is higher for all fixed $\alpha$ scenarios~(\autoref{tab:soedpe-stats}).

In summary, fixing $\alpha$ leads to a very constrained experimental design, and, if $\alpha$ is fixed to wrong values - the resulting model is biased.

\subsubsection{Regularization}

The results below are generated using the sOED-PE method with and without regularization of the PE and OED optimization problems.
In~\autoref{tab:soedpe-reg-stats} different metrics are summarized for all mixtures and the default precision accuracy scenario $\text{SCE}_{\sigma_y}$~(\autoref{tab:meas-quality}) for 6 different regularization scenarios: without regularization, with 4 different regularization techniques applied only to PE, with \emph{GO} applied both to PE and OED.
The change of the prediction quality, accuracy, and number of identifiable parameters over the sOED-PE iterations for one selected scenario is given in~\autoref{fig:regularization-diff}.

\begin{table}
    \centering
    \caption{Average and standard deviation of number of identifiable parameters $\overline{N}_i$ and prediction quality $Q^{0.95}_{\sigma_{x^V_1}}$, and average prediction accuracy $\overline{\sigma}_{x_1^V}$ for all mixtures for the first and last sOED-PE iteration $I_{sOED-PE}$. Discrepancy $D^{U}$ is calculated for all experiments. Different parameter regularization  $\text{SCE}_{Reg}$ scenario. In scenario $\text{OED}^{OED}$ PE and OED were regularized, in $\text{OED}^{OED}$ no regularization was used. In all other regularization scenario only PE was regularized. Default measurement accuracy scenario $\text{SCE}_{\sigma_y}$.}
    \label{tab:soedpe-reg-stats}
\begin{tabular}{cc|ccccc}
\toprule
\rot{$\text{SCE}_{Reg}$}           & \rot{$I_{sOED-PE}$} & $\overline{N}_{I}$ & $Q^{0.95}_{x_1^V}$ & \rot{$\overline{\sigma}_{x_1^V}$, \texttimes$10^{-4}$} & \rot{$D^{U}$, \texttimes$10^{-2}$} \\
\midrule                                                                      
\multirow{3}{*}{No Reg}            & 1                   & -                  & 1.00 $\pm$ 0.02    & 5.5                                                    & - \\
                                   & 15                  & -                  & 0.97 $\pm$ 0.09    & 2.0                                                    & - \\
                                   & All                 & -                  & -                  & -                                                      & 3.1\\
\cline{1-6}                                                                   
\multirow{3}{*}{E}                 & 1                   & 4.87 $\pm$ 0.34    & 1.00 $\pm$ 0.02    & 5.3                                                    & - \\
                                   & 15                  & 4.82 $\pm$ 0.39    & 0.96 $\pm$ 0.10    & 2.0                                                    & - \\
                                   & All                 & -                  & -                  & -                                                      & 2.9\\
\cline{1-6}                                                                   
\multirow{3}{*}{GO}                & 1                   & 3.90 $\pm$ 0.69    & 0.98 $\pm$ 0.06    & 4.2                                                    & - \\
                                   & 15                  & 4.08 $\pm$ 0.75    & 0.89 $\pm$ 0.17    & 1.8                                                    & - \\
                                   & All                 & -                  & -                  & -                                                      & 2.6\\
\cline{1-6}                                                                   
\multirow{3}{*}{$\text{GO}^{OED}$} & 1                   & 3.89 $\pm$ 0.69    & 0.98 $\pm$ 0.06    & 4.2                                                    & - \\
                                   & 15                  & 3.94 $\pm$ 0.64    & 0.87 $\pm$ 0.19    & 1.7                                                    & - \\
                                   & All                 & -                  & -                  & -                                                      & 5.9\\
\cline{1-6}                                                                   
\multirow{3}{*}{SVD}               & 1                   & 3.81 $\pm$ 0.71    & 0.97 $\pm$ 0.09    & 4.1                                                    & - \\
                                   & 15                  & 3.62 $\pm$ 0.74    & 0.86 $\pm$ 0.20    & 1.7                                                    & - \\
                                   & All                 & -                  & -                  & -                                                      & 2.4\\
\cline{1-6}                                                                   
\multirow{3}{*}{FS}                & 1                   & 2.41 $\pm$ 0.75    & 0.98 $\pm$ 0.08    & 4.2                                                    & - \\
                                   & 15                  & 2.79 $\pm$ 0.63    & 0.89 $\pm$ 0.18    & 1.8                                                    & - \\
                                   & All                 & -                  & -                  & -                                                      & 2.7\\
\cline{1-6}
\bottomrule
\end{tabular}
\end{table}

The difference in performance of the regularization techniques discussed for the first case study was observed again: based on the average number of identifiable parameters $\overline{N}_i$ the techniques are ranked in the same order: \emph{E}, \emph{GO}, \emph{SVD}, and \emph{FS}.
However, based on the standard deviation of the prediction quality $Q^{0.95}_y$, even the least conservative technique \emph{E} failed to consistently regularize without loss of prediction quality.
We assume that this is mainly due to the fact that OED with A criterion selected the experimental design that resulted in ``overfitted'' models after the final sOED-PE iteration.
This can be seen in the high values of the standard deviation ($\pm0.09$~\autoref{tab:soedpe-reg-stats}) of the prediction quality $Q^{0.95}_{x_1^V}$ at the final sOED-PE iteration even without the use of regularization $\left(\text{SCE}_{Reg}=\text{No Reg}\right)$.

There is no clear effect of the regularization techniques on the discrepancy of the experimental space as seen in~\autoref{tab:soedpe-reg-stats}.
The only apparent effect is due to regularization of OED $\left(\text{SCE}_{Reg}=\text{GO}^{OED}\right)$: the experimental design space becomes smaller, as expected from the analysis of the consequences of fixing $\alpha$.
The performance of the \emph{GO} regularization with $\left(\text{SCE}_{Reg}=\text{GO}^{OED}\right)$ and without regularization of OED $\left(\text{SCE}_{Reg}=\text{GO}\right)$ is shown for measurement accuracy scenarios $\text{SCE}_{\sigma_y}$ ``default''~\autoref{tab:soed-all-default} and ``precise''~\autoref{tab:soed-all-precise} .
The scenario $\text{SCE}_{Reg}=\text{GO}^{OED}$ performs similarly to the scenario $\text{SCE}_{Reg}=\text{GO}$, despite the increased discrepancy $D^{U}$, i.e., smaller experimental design space, suggesting that regularization of the OED makes the optimization problem easier to solve without introducing significant bias.

In general, sOED-PE performed well with artificial data and without the use of regularization techniques, except for the regularization used by the IPOPT solver itself, which was also seen in a real VLE experiment~\citep{Bubel2024}.
However, if the regularization technique must be used to overcome identifiability issues, it is recommended that either \emph{E} or \emph{GO} techniques be used, as they introduce the least regularization bias on average.

\section{Conclusions}

In this work, we have studied the identifiability of the NRTL parameters in binary vapor-liquid equilibria to demonstrate the shortcomings of the established heuristic of fixing $\alpha$, and to propose better alternatives.
For this purpose, existing methods for the analysis and computation of the exact model uncertainty based on Monte Carlo methods were modified to include the influence of regularization techniques.
The methods were then applied to several binary mixtures with different azeotropic behaviors.
In general, all 5 NRTL parameters could always be estimated independently of the binary mixture and measurement accuracy.
In addition, it was shown that the NRTL model can be highly nonlinear in parameters depending on the selected mixture and measurement accuracy, i.e., the real parameter uncertainty does not fit into the elliptical confidence region approximated by the linearization assumption.
However, despite the parameter nonlinearity, the real prediction accuracy of the model was correctly approximated by the linearization assumption.

In the first case study, it was shown that fixing $\alpha$ to the wrong value as the result of a heuristic can lead to a biased model, i.e., the model predictions are incorrect.
In some cases, the resulting model bias may be difficult to detect.
Of course, when $\alpha$ is fixed, the parameter estimation problem is easier to solve: there are only 4 parameters left to estimate, and the model becomes less nonlinear, though sometimes still ill-conditioned (overstretched covariance ellipses).
However, using the modern optimization solver IPOPT with internal regularization for PE, we always found the local optima, so estimating all 5 NRTL parameters simultaneously should be possible with modern solvers, provided the initial guess for the optimizer is close enough to the solution.
In summary, we strongly advise against fixing $\alpha$ due to the risks of a biased model, and emphasize, that additional experimental data must be used to validate the model if $\alpha$ is fixed to detect the bias.

The only alternative to fixing $\alpha$ is to estimate all 5 NRTL parameters simultaneously.
Estimation of all parameters may require the use of regularization, which improves prediction accuracy but introduces bias, as any regularization technique does.
Based on our studies, we recommend using an optimization solver with internal regularization for the initial parameter estimation and the Generalized Orthogonalization (\emph{GO}) regularization technique to identify the set of identifiable parameters, as it is independent of a user-defined threshold and has proven to be the most effective technique, balancing prediction accuracy and prediction bias.

The second case study highlighted the importance and applicability of the sOED-PE method to binary VLE measurements, even without regularization, as already shown in other studies~\citep{Duarte2021, Bubel2024}.
In general, the use of the A criterion used for sOED-PE resulted in a constant decrease of the prediction quality.
We assume that this is because many proposed experiments were concentrated around certain design regions, and as the number of experiments increased, the model became ``overfitted'' in these regions.
This behavior should persist for other classical OED criteria based only on the parameter variance--covariance matrix, e.g., D, E, etc.
Overall, it was shown that fixing $\alpha$, even to correct values, leads to an even larger reduction of the experimental design space, i.e., only a few unique experiments are suggested to be executed.
This, together with the risk of obtaining a biased model by fixing $\alpha$ to a wrong value, suggests that $\alpha$ should not be fixed when used together with sOED-PE.
The sOED-PE performed best when the internal regularization of the NLP solver was used without any additional regularization technique.
If a regularization technique must be used, i.e., the underlying PE or OED optimization problem cannot be solved due to identifiability problems, either \emph{E} or \emph{GO} techniques should be used, as they introduce the least regularization bias on average.
It was also shown that the experimenter can directly influence the discrepancy of the experimental space, i.e. the proposed experiments are closer to each other, by using regularization techniques with OED.

To further investigate this issue, the existing sOED-PE analysis methods need to include a globalization strategy for PE and different OED criteria that aim to improve not only the parameter accuracy but also the prediction accuracy.
In addition, although our results were consistent over 12 different mixtures with different azeotropic behaviors, it would be interesting to include more mixtures in the analysis to be more confident about the performance of different regularization techniques.
Finally, since the NRTL model is also used to describe multicomponent and liquid-liquid equilibrium, it would be interesting to apply the methods proposed in this work to other cases to see if the identifiability of the model worsens and the use of regularization techniques becomes the only feasible alternative to estimating all NRTL parameters.

\section{Funding Acknowledgement}
Gefördert durch die Deutsche Forschungsgemeinschaft (DFG) 466397921. This work is funded by the Deutsche Forschungsgemeinschaft (DFG, German Research Foundation) SPP 2331 - 466397921.

\bibliographystyle{elsarticle-harv}
\bibliography{bib_used}

\appendix

\section{Appendix}

\clearpage

\begin{table*}
\centering
\caption{Summary of binary mixtures, source of NRTL parameters and mixture properties. Azeotrope types as defined by \citet{Gmehling2019}:  Type 0: no azeotrope; Type I: homogeneous azeotrope, pressure maximum; Type III: homogeneous azeotrope, pressure minimum; Type V: double azeotrope. NRTL and Wagner25 parameters are taken from~\citet{AspenTechnology}}
\label{tab:mixtures}
\begin{tabular}{lcccc}
\toprule
 & Source & $\alpha^\ast$ & Azeotrope Type & $A_{12}=A_{21}=0$ \\
Mixture components &  &  &  &  \\
\midrule
Acetone / Chloroform & APV121 VLE-IG & 0.30 & Type III &  \\
Methanol / Benzene & APV121 VLE-IG & 0.40 & Type I &  \\
Methanol / Water & APV121 VLE-IG & 0.30 & Type 0 &  \\
Benzene / Hexafluorobenzene & APV121 VLE-IG & 0.30 & Type V & Yes \\
Chloroform / Ethyl acetate & APV121 VLE-IG & 0.30 & Type III & Yes \\
Ethanol / Benzene & APV121 VLE-IG & 0.30 & Type I &  \\
Ethanol / Dioxane & APV121 VLE-IG & 0.30 & Type I &  \\
Acetone / Ethanol & APV121 VLE-IG & 0.30 & Type 0 &  \\
Methanol / Ethanol & APV121 VLE-IG & 0.30 & Type 0 &  \\
Methanol / Epoxybutane & NISTV121 NIST-IG & 0.10 & Type I &  \\
Chloroform / Epoxybutane & NISTV121 NIST-IG & 0.50 & Type III &  \\
Ethanol / Acetic acid & NISTV121 NIST-IG & 0.50 & Type 0 &  \\
\bottomrule
\end{tabular}
\end{table*}

\begin{table}
    \centering
    \caption{Average and standard deviation of prediction quality $Q^{0.95}_y$ and number of identifiable parameters $N_i$ for all mixtures for 4 regularization techniques $\text{SCE}_{Reg}$. Different scenarios for measurement variables $\text{SCE}_{y}$, and measurement accuracies $\text{SCE}_{\sigma_y}=[\sigma_{x^V_1},~\sigma_{T}]$. Precise: $[0.0002~\si{\mole\per\mole},~0.01~\si{K}]$. Default: $[0.001~\si{\mole\per\mole},~0.03~\si{K}]$.}

\label{tab:reg-statistics}
\begin{tabular}{ccccc}
\toprule
 &  &  & $Q^{0.95}_y$ & $N_i$ \\
$\text{SCE}_{y}$ & $\text{SCE}_{\sigma_y}$ & $\text{SCE}_{Reg}$ &  &  \\
\midrule
\multirow{8}{*}{\rot{$[x^V_1]$}} & \multirow{4}{*}{\rot{default}} & E & 0.95 $\pm$ 0.04 & 4.70 $\pm$ 0.61 \\
 &  & GO & 0.94 $\pm$ 0.03 & 4.17 $\pm$ 0.70 \\
 &  & SVD & 0.93 $\pm$ 0.03 & 4.03 $\pm$ 0.68 \\
 &  & FS & 0.88 $\pm$ 0.04 & 2.94 $\pm$ 0.84 \\
\cline{2-5}
 & \multirow{4}{*}{\rot{precise}} & E & 0.95 $\pm$ 0.04 & 4.83 $\pm$ 0.48 \\
 &  & GO & 0.94 $\pm$ 0.03 & 4.58 $\pm$ 0.56 \\
 &  & SVD & 0.91 $\pm$ 0.04 & 3.95 $\pm$ 0.80 \\
 &  & FS & 0.87 $\pm$ 0.04 & 2.99 $\pm$ 0.82 \\
\cline{1-5} \cline{2-5}
\multirow{8}{*}{\rot{$[x^V_1,T]$}} & \multirow{4}{*}{\rot{default}} & E & 0.94 $\pm$ 0.05 & 4.71 $\pm$ 0.74 \\
 &  & GO & 0.93 $\pm$ 0.05 & 4.37 $\pm$ 0.79 \\
 &  & SVD & 0.90 $\pm$ 0.05 & 3.86 $\pm$ 0.81 \\
 &  & FS & 0.86 $\pm$ 0.05 & 2.93 $\pm$ 0.87 \\
\cline{2-5}
 & \multirow{4}{*}{\rot{precise}} & E & 0.95 $\pm$ 0.03 & 4.82 $\pm$ 0.56 \\
 &  & GO & 0.94 $\pm$ 0.03 & 4.60 $\pm$ 0.60 \\
 &  & SVD & 0.91 $\pm$ 0.04 & 3.88 $\pm$ 0.89 \\
 &  & FS & 0.87 $\pm$ 0.04 & 2.98 $\pm$ 0.85 \\
\cline{1-5} \cline{2-5}
\bottomrule
\end{tabular}
\end{table}

\begin{table}
\centering
\caption{Average frequency with which parameters are marked as identifiable, where 1 means the parameter is always identifiable.}
\label{tab:ident-parameter}
\begin{tabular}{c|cccc}
\toprule
& \multicolumn{4}{c}{$\text{SCE}_{Reg}$} \\ 
     & E & GO & SVD & FS \\
\midrule
$A_{12}$ & 0.95 $\pm$ 0.12 & 0.86 $\pm$ 0.19 & 0.77 $\pm$ 0.34 & 0.69 $\pm$ 0.41 \\
$A_{21}$ & 0.95 $\pm$ 0.12 & 0.83 $\pm$ 0.23 & 0.76 $\pm$ 0.31 & 0.64 $\pm$ 0.40 \\
$B_{12}$ & 0.95 $\pm$ 0.12 & 0.92 $\pm$ 0.18 & 0.77 $\pm$ 0.35 & 0.64 $\pm$ 0.40 \\
$B_{21}$ & 0.95 $\pm$ 0.12 & 0.88 $\pm$ 0.21 & 0.76 $\pm$ 0.35 & 0.70 $\pm$ 0.43 \\
$\alpha$ & 0.96 $\pm$ 0.12 & 0.94 $\pm$ 0.14 & 0.86 $\pm$ 0.24 & 0.29 $\pm$ 0.37 \\
\bottomrule
\end{tabular}
\end{table}

\begin{table}
\centering
\caption{Number of scenarios (SCE) for every mixture used in \emph{Case Study II}. Special case $^a$: Methanol / Epoxybutane has only 16 runs because $\underline{\alpha} = \alpha^\ast$.}
\label{tab:cs2-scenarios}
\begin{tabular}{ccccc}
\toprule
 & $\text{SCE}_{\sigma}$ & $\text{SCE}_{Reg}$ & Total \\
$\text{SCE}_{\theta}$ &  &  &  \\
\midrule
\emph{All} & 2 & 6 & 12 \\
$\alpha^\ast$ & 2 & - & 2 \\
$\overline{\alpha}$ & 2 & - & 2 \\
$\underline{\alpha}$ & 2 & - & $2^a$ \\
\midrule
Total &  &  & 18 \\
\bottomrule
\end{tabular}
\end{table}

\begin{figure}[h]
    \centering
    \includegraphics[width=\columnwidth]{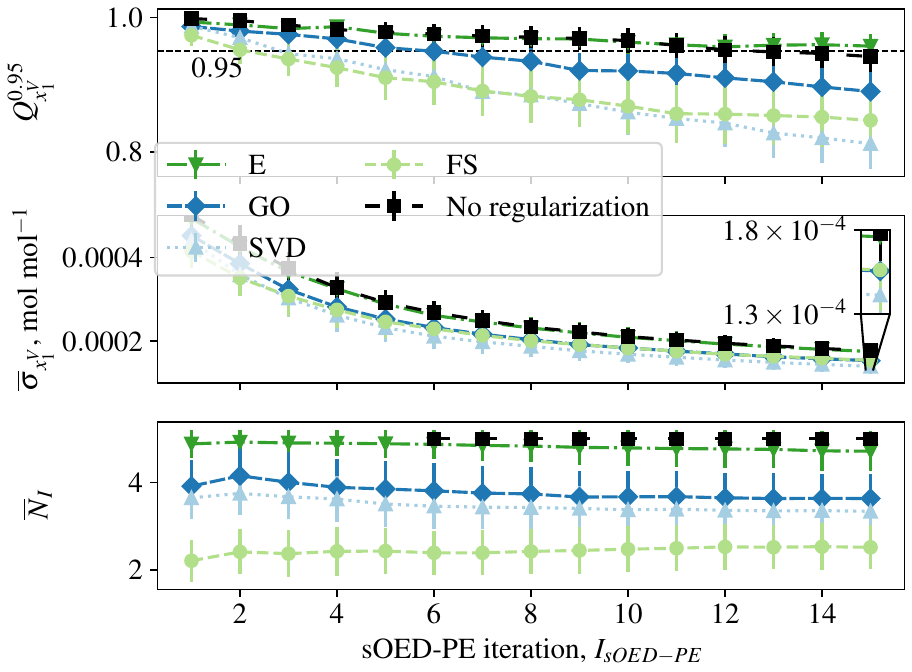}
    \caption{Change of prediction quality $Q^{0.95}_{x^V_1}$, accuracy $\overline{\sigma}^{0.95}_{x^V_1}$, and average number of identifiable parameters $\overline{N}_i$ over sOED-PE iterations $I_{sOED-PE}$ for different regularization scenarios. Chloroform / Etyalacetate mixture, both measurement variables are used $\text{SCE}_{y}$ with default measurement accuracy $\text{SCE}_{\sigma_y}$. Standard deviation of the prediction quality $Q^{0.95}_y$ is scaled down by a factor of 0.2 for better visualization.}
    \label{fig:regularization-diff}
\end{figure}

\begin{table*}
    \centering
    \caption{Average and standard deviation of prediction quality $Q^{0.95}_{\sigma_{x^V_1}}$, number of identifiable parameters $\overline{N}_i$ and average estimation accuracy  $\overline{s}_{x_1^V}$ and prediction accuracy $\overline{\sigma}_{x_1^V}$ for all mixtures for the first and last sOED-PE iteration $I_{sOED-PE}$. Discrepancy of the final design space for all mixtures $D^{U}$. Different parameter $\text{SCE}_{\theta}$ and regularization  $\text{SCE}_{Reg}$ scenario. In scenario $\text{OED}^{OED}$ PE and OED were regularized. In all other regularization scenario only PE was regularized. Precise measurement accuracies: $\sigma^M_{x^V_1}=0.0002~\si{\mole\per\mole},\sigma^M_{T}=0.01~\si{K}]$.}
    \label{tab:soed-all-precise}
\begin{tabular}{cccccccc}
\toprule
 &  &  & $Q^{0.95}_{x_1^V}$ & $\overline{s}_{x_1^V}$ & $\overline{\sigma}_{x_1^V}$ & $\overline{N}_{I}$ & $D^{U}$ \\
\rot{$\text{SCE}_{\theta}$} & \rot{$\text{SCE}_{Reg}$} & \rot{$I_{sOED-PE}$} &  &  &  &  &  \\
\midrule
\multirow{3}{*}{$\alpha^\ast$} & \multirow{3}{*}{-}  & 1& 0.99 $\pm$ 0.04 & 1.8\texttimes$10^{-4}$ & 9.8\texttimes$10^{-5}$ & - & - \\
 &  & 15 & 0.96 $\pm$ 0.10 & 1.9\texttimes$10^{-4}$ & 4.2\texttimes$10^{-5}$ & - & - \\
 &  & All & - & - & - & - & 9.0\texttimes$10^{-2}$ \\
\cline{1-8} \cline{2-8}
\multirow{3}{*}{$\underline{\alpha}$} & \multirow{3}{*}{-}  & 1& 0.84 $\pm$ 0.18 & 2.2\texttimes$10^{-4}$ & 1.0\texttimes$10^{-4}$ & - & - \\
 &  & 15 & 0.48 $\pm$ 0.29 & 2.3\texttimes$10^{-4}$ & 4.1\texttimes$10^{-5}$ & - & - \\
 &  & All & - & - & - & - & 6.9\texttimes$10^{-2}$ \\
\cline{1-8} \cline{2-8}
\multirow{3}{*}{$\overline{\alpha}$} & \multirow{3}{*}{-}  & 1& 0.78 $\pm$ 0.20 & 2.2\texttimes$10^{-4}$ & 1.1\texttimes$10^{-4}$ & - & - \\
 &  & 15 & 0.36 $\pm$ 0.25 & 2.3\texttimes$10^{-4}$ & 4.2\texttimes$10^{-5}$ & - & - \\
 &  & All & - & - & - & - & 9.0\texttimes$10^{-2}$ \\
\cline{1-8} \cline{2-8}
\multirow{18}{*}{\emph{All}} & \multirow{3}{*}{-}  & 1& 1.00 $\pm$ 0.02 & 1.7\texttimes$10^{-4}$ & 1.4\texttimes$10^{-4}$ & - & - \\
 &  & 15 & 0.97 $\pm$ 0.09 & 1.9\texttimes$10^{-4}$ & 5.4\texttimes$10^{-5}$ & - & - \\
 &  & All & - & - & - & - & 3.7\texttimes$10^{-2}$ \\
\cline{2-8}
 & \multirow{3}{*}{E}  & 1& 1.00 $\pm$ 0.03 & 1.7\texttimes$10^{-4}$ & 1.3\texttimes$10^{-4}$ & 4.95 $\pm$ 0.22 & - \\
 &  & 15 & 0.96 $\pm$ 0.10 & 1.9\texttimes$10^{-4}$ & 5.4\texttimes$10^{-5}$ & 4.91 $\pm$ 0.30 & - \\
 &  & All & - & - & - & - & 3.6\texttimes$10^{-2}$ \\
\cline{2-8}
 & \multirow{3}{*}{GO}  & 1& 0.98 $\pm$ 0.06 & 1.7\texttimes$10^{-4}$ & 1.2\texttimes$10^{-4}$ & 4.52 $\pm$ 0.61 & - \\
 &  & 15 & 0.93 $\pm$ 0.16 & 1.9\texttimes$10^{-4}$ & 5.2\texttimes$10^{-5}$ & 4.57 $\pm$ 0.65 & - \\
 &  & All & - & - & - & - & 3.4\texttimes$10^{-2}$ \\
\cline{2-8}
 & \multirow{3}{*}{$\text{GO}^{OED}$}  & 1& 0.98 $\pm$ 0.06 & 1.7\texttimes$10^{-4}$ & 1.2\texttimes$10^{-4}$ & 4.53 $\pm$ 0.61 & - \\
 &  & 15 & 0.92 $\pm$ 0.17 & 1.9\texttimes$10^{-4}$ & 5.1\texttimes$10^{-5}$ & 4.53 $\pm$ 0.61 & - \\
 &  & All & - & - & - & - & 5.7\texttimes$10^{-2}$ \\
\cline{2-8}
 & \multirow{3}{*}{SVD}  & 1& 0.95 $\pm$ 0.11 & 1.6\texttimes$10^{-4}$ & 1.0\texttimes$10^{-4}$ & 3.75 $\pm$ 0.69 & - \\
 &  & 15 & 0.81 $\pm$ 0.25 & 2.1\texttimes$10^{-4}$ & 4.7\texttimes$10^{-5}$ & 3.79 $\pm$ 0.79 & - \\
 &  & All & - & - & - & - & 3.5\texttimes$10^{-2}$ \\
\cline{2-8}
 & \multirow{3}{*}{FS}  & 1& 0.96 $\pm$ 0.10 & 1.6\texttimes$10^{-4}$ & 1.0\texttimes$10^{-4}$ & 2.36 $\pm$ 0.75 & - \\
 &  & 15 & 0.84 $\pm$ 0.22 & 2.0\texttimes$10^{-4}$ & 4.7\texttimes$10^{-5}$ & 2.72 $\pm$ 0.70 & - \\
 &  & All & - & - & - & - & 3.6\texttimes$10^{-2}$ \\
\cline{1-8} \cline{2-8}
\bottomrule
\end{tabular}
\end{table*}

\begin{table*}
    \centering
    \caption{Average and standard deviation of prediction quality $Q^{0.95}_{\sigma_{x^V_1}}$, number of identifiable parameters $\overline{N}_i$ and average estimation accuracy  $\overline{s}_{x_1^V}$ and prediction accuracy $\overline{\sigma}_{x_1^V}$ for all mixtures for the first and last sOED-PE iteration $I_{sOED-PE}$. Discrepancy of the final design space for all mixtures $D^{U}$. Different parameter $\text{SCE}_{\theta}$ and regularization  $\text{SCE}_{Reg}$ scenario. In scenario $\text{OED}^{OED}$ PE and OED were regularized. In all other regularization scenario only PE was regularized. Default measurement accuracies: $\sigma_{x^V_1}=0.001~\si{\mole\per\mole},\sigma_{T}=0.03~\si{K}]$.}
    \label{tab:soed-all-default}
\begin{tabular}{cccccccc}
\toprule
 &  &  & $Q^{0.95}_{x_1^V}$ & $\overline{s}_{x_1^V}$ & $\overline{\sigma}_{x_1^V}$ & $\overline{N}_{I}$ & $D^{U}$ \\
\rot{$\text{SCE}_{\theta}$} & \rot{$\text{SCE}_{Reg}$} & \rot{$I_{sOED-PE}$} &  &  &  &  &  \\
\midrule
\multirow{3}{*}{$\alpha^\ast$} & \multirow{3}{*}{-} & 1 & 0.99 $\pm$ 0.05 & 1.0\texttimes$10^{-3}$ & 3.8\texttimes$10^{-4}$ & - & - \\
 &  & 15 & 0.96 $\pm$ 0.11 & 9.9\texttimes$10^{-4}$ & 1.5\texttimes$10^{-4}$ & - & - \\
 &  & All & - & - & - & - & 6.7\texttimes$10^{-2}$ \\
\cline{1-8} \cline{2-8}
\multirow{3}{*}{$\underline{\alpha}$} & \multirow{3}{*}{-} & 1 & 0.94 $\pm$ 0.15 & 1.0\texttimes$10^{-3}$ & 4.0\texttimes$10^{-4}$ & - & - \\
 &  & 15 & 0.74 $\pm$ 0.28 & 1.0\texttimes$10^{-3}$ & 1.4\texttimes$10^{-4}$ & - & - \\
 &  & All & - & - & - & - & 4.9\texttimes$10^{-2}$ \\
\cline{1-8} \cline{2-8}
\multirow{3}{*}{$\overline{\alpha}$} & \multirow{3}{*}{-} & 1 & 0.93 $\pm$ 0.13 & 1.0\texttimes$10^{-3}$ & 4.0\texttimes$10^{-4}$ & - & - \\
 &  & 15 & 0.72 $\pm$ 0.25 & 1.0\texttimes$10^{-3}$ & 1.5\texttimes$10^{-4}$ & - & - \\
 &  & All & - & - & - & - & 6.6\texttimes$10^{-2}$ \\
\cline{1-8} \cline{2-8}
\multirow{18}{*}{\emph{All}} & \multirow{3}{*}{-} & 1 & 1.00 $\pm$ 0.02 & 1.0\texttimes$10^{-3}$ & 5.5\texttimes$10^{-4}$ & - & - \\
 &  & 15 & 0.97 $\pm$ 0.09 & 1.0\texttimes$10^{-3}$ & 2.0\texttimes$10^{-4}$ & - & - \\
 &  & All & - & - & - & - & 3.1\texttimes$10^{-2}$ \\
\cline{2-8}
 & \multirow{3}{*}{E} & 1 & 1.00 $\pm$ 0.02 & 1.0\texttimes$10^{-3}$ & 5.3\texttimes$10^{-4}$ & 4.87 $\pm$ 0.34 & - \\
 &  & 15 & 0.96 $\pm$ 0.10 & 1.0\texttimes$10^{-3}$ & 2.0\texttimes$10^{-4}$ & 4.82 $\pm$ 0.39 & - \\
 &  & All & - & - & - & - & 2.9\texttimes$10^{-2}$ \\
\cline{2-8}
 & \multirow{3}{*}{GO} & 1 & 0.98 $\pm$ 0.06 & 9.6\texttimes$10^{-4}$ & 4.2\texttimes$10^{-4}$ & 3.90 $\pm$ 0.69 & - \\
 &  & 15 & 0.89 $\pm$ 0.17 & 1.0\texttimes$10^{-3}$ & 1.8\texttimes$10^{-4}$ & 4.08 $\pm$ 0.75 & - \\
 &  & All & - & - & - & - & 2.6\texttimes$10^{-2}$ \\
\cline{2-8}
 & \multirow{3}{*}{$\text{GO}^{OED}$} & 1 & 0.98 $\pm$ 0.06 & 9.5\texttimes$10^{-4}$ & 4.2\texttimes$10^{-4}$ & 3.89 $\pm$ 0.69 & - \\
 &  & 15 & 0.87 $\pm$ 0.19 & 1.0\texttimes$10^{-3}$ & 1.7\texttimes$10^{-4}$ & 3.94 $\pm$ 0.64 & - \\
 &  & All & - & - & - & - & 5.9\texttimes$10^{-2}$ \\
\cline{2-8}
 & \multirow{3}{*}{SVD} & 1 & 0.97 $\pm$ 0.09 & 9.5\texttimes$10^{-4}$ & 4.1\texttimes$10^{-4}$ & 3.81 $\pm$ 0.71 & - \\
 &  & 15 & 0.86 $\pm$ 0.20 & 1.0\texttimes$10^{-3}$ & 1.7\texttimes$10^{-4}$ & 3.62 $\pm$ 0.74 & - \\
 &  & All & - & - & - & - & 2.4\texttimes$10^{-2}$ \\
\cline{2-8}
 & \multirow{3}{*}{FS} & 1 & 0.98 $\pm$ 0.08 & 9.6\texttimes$10^{-4}$ & 4.2\texttimes$10^{-4}$ & 2.41 $\pm$ 0.75 & - \\
 &  & 15 & 0.89 $\pm$ 0.18 & 1.0\texttimes$10^{-3}$ & 1.8\texttimes$10^{-4}$ & 2.79 $\pm$ 0.63 & - \\
 &  & All & - & - & - & - & 2.7\texttimes$10^{-2}$ \\
\cline{1-8} \cline{2-8}
\bottomrule
\end{tabular}
\end{table*}

\end{document}